\documentclass[journal,11pt,twocolumn]{IEEEtran}
\usepackage{graphicx}
\usepackage{amsmath}

\newcommand{\ignore}[1]{}
\usepackage{subfigure}
\usepackage{subeqnarray}
\usepackage{amsfonts}
\usepackage{mathrsfs}
\usepackage{lineno, hyperref}
\usepackage{textcomp}
\usepackage{xcolor}
\usepackage{multirow}
\usepackage{comment}
\graphicspath{ {figures/} }
\usepackage{float}
\usepackage{soul}

\usepackage[justification=centering]{caption}

\usepackage{cite}

 \usepackage{dblfloatfix}

\begin{document}
%
\title{Formation Mission Design for Commercial Aircraft using Switched Optimal Control Techniques}
%
%
%

\author{Mar\'{\i}a Cerezo-Maga\~na, 
        Alberto Olivares, 
        Ernesto Staffetti
\thanks{{Universidad Rey Juan Carlos, Department of Telecommunication Engineering, Fuenlabrada, Madrid, Spain (e-mail: maria.cerezo@urjc.es; alberto.olivares@urjc.es; ernesto.staffetti@urjc.es).}}
}

\markboth{IEEE Transactions on Aerospace and Electronic Systems}%
{Cerezo-Maga\~na \MakeLowercase{\textit{et al.}}: Formation Mission Design for Commercial Aircraft using Switched Optimal Control Techniques}

\maketitle

\begin{abstract}
In this paper, the formation mission design problem for commercial
aircraft is studied.
Given the departure times and the departure and arrival locations of
several commercial flights, the relevant weather forecast, and the
expected fuel savings during formation flight, the problem consists in
establishing how to organize them in formation or solo flights and in
finding the trajectories that minimize the overall direct operating cost
of the flights.
Since each aircraft can fly solo or in different positions inside a
formation, the mission is modeled as a switched dynamical system, in
which the discrete state describes the combination of flight modes of
the individual aircraft and logical constraints in disjunctive form
establish the switching logic among the discrete states of the system.
The formation mission design problem has been formulated as an optimal
control problem of a switched dynamical system and solved using an
embedding approach, which allows switching decision among discrete
states to be modeled without relying on binary variables.
The resulting problem is a classical optimal control problem which has been solved using a
knotting pseudospectral method. Several numerical experiments have been
conducted to demonstrate the effectiveness of this approach. 
The obtained results show that formation flight has great potential to
reduce fuel consumption and emissions.
\end{abstract}


\IEEEpeerreviewmaketitle

\section{Introduction}

The aviation sector represents a great challenge for the environment, being a top-ten {global GreenHouse Gases} (\textsf{GHG}) emitter with a 4.9$\%$  contribution to the Earth\textquotesingle s warming \cite{lee2009aviation}.
According to a recent special report of the Intergovernmental Panel on Climate Change (\textsf{IPCC})  \cite{IPCC2018}, a leading climate science body, aviation is a crucial sector to the threat of climate change and global warming. 
Indeed, although air transport is \textcolor{black}{only} responsible for 14\% of the dioxide carbon emissions of the transport sector, much less than other transport modes like light-duty vehicles (30\%) and heavy-duty vehicles (36\%), it seems to be harder to decarbonized. 
The reason is the air traffic predicted growth rate, which is projected to be greater than other transport modes.
Indeed, the air traffic predicted growth rate is another concerning aspect. The Compound Annual Growth Rate (\textsf{CAGR}) forecast over the next two decades, assessed by the International Air Transport Association (\textsf{IATA}), estimates that the number of flights is expected to continue increasing by the 3.5$\%$ annually, which means doubling the current number of passengers  by 2040 \cite{international2016iata}. This continuous growth in air traffic demand imposes the search for new solutions to mitigate its potential negative effects on our planet.
Unfortunately, 
electrification of large aircraft is unlikely to happen anytime soon due to problems related to power to weight ratio and power density levels of the batteries
 \cite{madonna2018electrical}.

Hence, the predicted effects of the aviation sector to climate change are greater than desired \cite{IPCC2018}, \cite{icao2016environmental}, and a great effort is required to mitigate them. In this sense, formation flight has the potential to make a significant contribution.
Even more, formation flight offers a great promise 
not only in mitigating the environmental impact of aviation, but also, in increasing the {Air Traffic Management}  (\textsf{ATM}) capacity.
According to \cite{chicagoconvention}, ``\textit{the formation operates as a single aircraft with regard to navigation and position reporting}'', which means that two or three aircraft flying in formation would be treated as one single aircraft for air traffic control purposes. This would help to reduce their workload. 

\textcolor{black}{
The considered formation mission design problem can be stated as follows: given the departure times and the departure and arrival locations of several commercial flights, the relevant weather forecast, and the expected fuel savings during formation flight, the problem consists in establishing how to organize them in formation or solo flights and in finding the trajectories that minimize the overall Direct Operating Costs (\textsf{DOC}) of the flights. The rendezvous and splitting locations and times must also be determined.
}
\textcolor{black}{
\textsf{DOC} are all costs incurred directly in the operation of the aircraft. These expenses are usually divided into three groups: flight operating costs, such as flight crew salaries, fuel and oil, and airport and en-route charges, costs related to maintenance and overhaul, and aircraft depreciation \cite{camilleri:2018:tmteatap}. 
}

\textcolor{black}{
Only formations of up to three aircraft have been considered in this paper. This choice is motivated by the fact that the difficulties in synchronizing more than three flights make formations of more than three aircraft operationally impractical. Besides, increasing the number of aircraft in the formation reduces asymptotically the benefits obtained from the formation \cite{ningetal:2011:apoeff}.}

\textcolor{black}{
In this paper, the formation mission design problem has been formulated as an optimal control problem of a switched dynamical system. These systems are described by both a continuous and a discrete dynamics in which the transitions among discrete states are not established in advance. In particular, each aircraft has different flight modes, namely solo flight and flight in different positions inside a formation, and their combination is represented by the discrete state of the switched dynamical system, which models their joint dynamic behavior. Each flight mode will be represented by different dynamical equations which may include or not formation flight benefits in terms of fuel savings. Additionally, logical constraints in disjunctive form, based on the stream-wise distance between aircraft, model the switching logic among the discrete states of the system.}

\subsection{{Previous Approaches}}

The notion of formation flight arose from the observation of nature, in which many species of bird fly in formation during their migrations \cite{lissamanandshollenberger:1970:ffob}. 
Later, it was demonstrated that formation flight led to notable efficiency improvements for the birds and that substantial fuel savings were also possible for aircraft flying in formation. There are three types of formations: the echelon or in-line formation, the V formation and the inverted V formation.
These first studies mainly focused on the aerodynamic of formation flight to assess the potential benefits in terms of induced drag reductions and, consequently, in potential fuel burn reductions
\cite{ningetal:2011:apoeff,  kless2013inviscid, slotnick2014computational}.
Experimental studies on formation flight have also been carried out. One of the most exhaustive works has been conducted within the Surfing Aircraft Vortices for Energy (\textsf{\$AVE}) project \cite{bieniawskietal:2014:softarftfffabp}. The results obtained using two C-17s aircraft, at longitudinal distances of between 18 and 70 wingspans, have shown fuel savings of up to  11\% \cite{haalasetal:2014:fffabsdav, flanzeretal:2014:oaftfffabp} for the trailing aircraft. The results of another experimental study \cite{vachon2002f}, obtained using two F/A-18s flying closer, at a longitudinal distance of up to 6.6 wingspans, have shown fuel savings of just over 18\% for the trailing aircraft.

Besides these studies on aerodynamics of formation flight focused on
estimating the potential induced drag and fuel burn reductions, 
other research activities have been conducted to study different practical aspects of formation flight, such as the formation tracking control problem under actuator and sensor faults \cite{liu2019integrated} and the collision avoidance in formation flight \cite{seo2017collision}.

Regarding the optimal positioning, the lateral and the
stream-wise distance gain significance, in particular, the longitudinal
spacing between aircraft, usually expressed in terms of wingspans, is
used to classify formations. 
In the literature, the formations in which the longitudinal distance
between leader and trailing aircraft is less than 10 wingspans are
called close formations whereas those in which this distance is between
10 and 40 wingspans are called extended formations
\cite{durangoetal:2016ffifhefcta}. Close formations provide the greatest
benefits in terms of induced drag reduction while, in extended
formations, the potential induced drag reduction is lower, especially in
separations of more than 20 wingspans. However, extended formations are
intrinsically safer than close formations, a crucial aspect in
commercial flights.

Due to its potential in reducing fuel consumption and \textsf{GHG}
emissions, formation flight has recently attracted increasing interest
and the studies on aerodynamics of formation flight have given way to studies on formation flight planning \cite{xuetal:2014:arofff, kent:2015:aatorfcff, blakeandflanzer:2016:orfdrffarc, boweretal:2009:fgarofcff, hartjesetal:2018:mptofffica, hartjes2019trajectory}, most of which focus on extended formations. 
In these works, the formation flight planning problem has been solved as a bi-level optimization problem, in which the higher level problem is the optimal partner allocation problem and the lower level problem is the mission design problem.
While the former problem consists in establishing how to group the elements of a large-size set of flights into smaller subsets which contain potential partners of beneficial formation missions in terms of overall \textsf{DOC}, 
the latter problem consists in determining the optimal trajectories of the aircraft that compose candidate formation missions.

Although the mission design and the partner allocation problems are intrinsically very different, due to the existing interconnection between them, they have been addressed in one go in several works such as \cite{xuetal:2014:arofff} and \cite{kent:2015:aatorfcff}, in which the Breguet range equation and the Fermat point extension problem have been employed, respectively. Other examples are \cite{blakeandflanzer:2016:orfdrffarc, boweretal:2009:fgarofcff}, in which geometric methods have been used. 
In these approaches, the inherent combinatorial complexity of the partner allocation problem, which quickly grows with the number of flights, forced to introduce simplifications, 
such as using approximated dynamic models of the aircraft and neglecting the meteorological forecast, which are essential factors to achieve a satisfactory predictability of the trajectories and realistic estimations of the fuel savings. The natural framework to take into consideration all these elements is optimal control.
In this work only the mission design problem will be addressed using optimal control techniques and accurate dynamic models of the aircraft. The partner allocation problem is left for future research.

\textcolor{black}{Even though it is possible, choosing the type of formation and the relative position of each aircraft in the formation has not been included in the formulation of the formation mission design problem. Thus, without loss of generality, the echelon formation has been assumed and the position of each aircraft in the formation has been established in advance.
}

In \cite{hartjesetal:2018:mptofffica} the mission design problem has been studied for two-aircraft formations using a multi-phase optimal control approach. In this case, the number of phases and the structure of the formation missions is unique and known in advance. There are five phases. The first two correspond to the flights from the departure locations to the rendezvous point, after which the formation flight, which corresponds to the third phase, starts. The other two flight phases correspond to the flights from the splitting point to the arrival locations. In this case, a solution is found using multi-phase optimal control techniques and compared with the solution of a solo flight mission to establish which gives rise to minimum \textsf{DOC}. 
\textcolor{black}{
The same framework used in \cite{hartjesetal:2018:mptofffica} is extended in \cite{hartjes2019trajectory} to design three-aircraft formation missions. 
However, only in the latter a wind field model, obtained from a weather database using a polynomial regression, has been taken into account.}
In contrast to two-aircraft formation missions, the structure of three-aircraft formation missions is not unique and, consequently, many different cases must be considered, including the case in which one or all the aircraft flight solo. The obtained solutions must be then compared to establish which corresponds to the minimum \textsf{DOC}. This drawback is due to the fact that the multi-phase optimal control approach is only able to tackle optimal control problems with a fixed switching structure, that is, with transitions among phases established in advance. 

This approach is feasible if the number of flights considered is low 
and the type of formation and the relative position of each aircraft in the formation are fixed. 
Otherwise, the number of combinations quickly grows with the number of flights which makes the previous approach hardly scalable.

A more general framework is needed, which is able to design formation missions solving a single optimal control problem in which neither the number of phases nor the switching structure  \textcolor{black}{are} known in advance. Moreover, this approach should be easily scalable to design formation missions involving more than three flights and should allow a more accurate aerodynamic model to be used.

\subsection{{Contributions of the paper}}
In this paper, an innovative approach to the formation mission design problem is proposed which is formulated
as an optimal problem problem of a switched dynamical system with logical constraints in disjunctive form and solved using an embedding approach \cite{bengeaanddecarlo:2005:ocoss, bengeaetal:2011:ocossveicocp}.

The optimal control formulation of the problem allows accurate dynamic models of aircraft and meteorological forecast to be included in the problem formulation in order to improve the predictability of the trajectories and the estimation of the fuel savings. This represents a first improvement with respect to previous approaches \cite{xuetal:2014:arofff}, \cite{kent:2015:aatorfcff}, \cite{blakeandflanzer:2016:orfdrffarc}, \cite{boweretal:2009:fgarofcff}, and \cite{hartjesetal:2018:mptofffica}. Indeed, in the latter, which is based on optimal control, the meteorological forecast has not been taken into account.

Modeling the formation mission as a switched dynamical systems allows the multiphase approach adopted in 
\cite{hartjesetal:2018:mptofffica} and \cite{hartjes2019trajectory} to be avoided. This represents another  improvement of the method presented in this paper with respect to these previous approaches, since the embedding method employed to solve the resulting optimal control problem of a switched dynamical system
is able to deal with both the switching dynamics of the system and the logical constraints in disjunctive form in an efficient way. 
Additionally, the switching logic among discrete states can be modeled without using binary variables, 
the number of switches does not have to be established in advance, 
and the optimal values of the switching times between discrete states are obtained without introducing them as unknowns of the optimal control problem.  Therefore, the resulting problem is a classical optimal control problem, which has been solved using a pseudospectral knotting method \cite{Ross2004PseudospectralKM, Garg2017AnOO}.

The advantages in terms of computational time to solve the formation mission design problem given by the method presented in this paper with respect the multiphase approach have been quantified comparing the computational times obtained solving the same instances of the problem with both techniques.

The embedding approach has been used in \cite{seenivasan2019multiaircraft} to tackle logical constraints 
in disjunctive form introduced to enforce separation among aircraft in a multiple aircraft trajectory planning problem. 
However, in  \cite{seenivasan2019multiaircraft} there are no switches in the 
dynamic model of the system, whereas, as already mentioned, in the formation mission design problem studied in this paper, the embedding approach is used to handle both, the switched dynamical system and the logical constraints in disjunctive form.

Several numerical simulations have been carried out to show the effectiveness of the proposed technique. In the
first numerical experiment, a two-aircraft transoceanic mission design problem has been solved, in which the departure time of one flight has been left free. The benefits of formation flight with respect to solo flights have also been quantified. In a second numerical experiment,   an analysis of the solution obtained in the previous experiment has been carried out to determine how delays in the departure times of each flight affect the formation flight. In the third numerical experiment, a three-aircraft transoceanic mission design problem has been solved. An analysis of the obtained solution has also been performed to determine how the fuel savings scheme affects the formation flight.

\begin{figure}[ht!]
\centering
\renewcommand{\figurename}{Fig.}
\includegraphics[width = 25em]{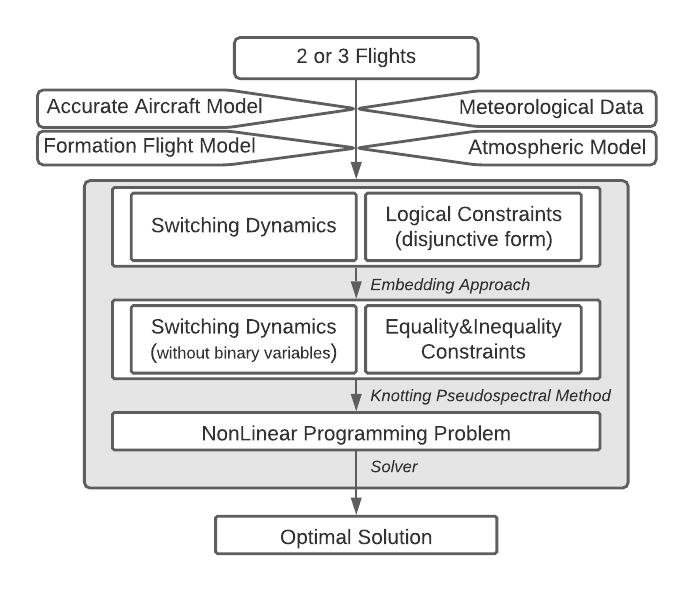}
\caption{{Schematic diagram of the formation mission design system.}}
\label{fig:system_architecture}
\end{figure}

The method for formation mission design presented in the paper is intended 
to be incorporated in a ground-based system {whose schematic diagram is given in Fig.~\ref{fig:system_architecture}, which could be used by flight dispatchers as a support system
 to design} formation missions for airlines or airline alliances.
Flight dispatchers are employed by airlines that are responsible for creating 
and modifying the flight plan. They also monitor the aircraft during flight. 
The computational time of the proposed algorithm makes it suitable for the strategical and tactical planning phases 
of the flight, as well as for the operational phase, that is, during flight.

\subsection{{Organization of the paper}}

The paper is organized as follows. In Section \ref{sect:model_of_the_system} the model of the aircraft is introduced \textcolor{black}{and} the modeling of the wind forecast is described. In the same section, the peculiarities of the formation mission design problem are discussed.
In Section \ref{sect:the_switched_optimal_control_problem}, the switched optimal control problem used to solve the mission design problem is stated. 
In Section \ref{sect:the_embedding_approach}, the embedding approach employed to solve the resulting switched optimal control problem is described. 
In Section \ref{sect:treatment_of_the_logical_constraints}, the technique employed to model logical constraints in disjunctive form is presented and then
particularized to model the logical constraints arising in the mission design problem. 
In Section \ref{sect:solution_of_the_switched_optimal_control_problem}, the pseudospectral knotting method used to solve the resulting classical optimal control problem is
outlined.
The results of the numerical experiments are reported and analyzed in Section \ref{sect:numerical_results}. 
Finally, in Section \ref{sect:conclusions}, some conclusions are drawn.

\section{Model of the system}
\label{sect:model_of_the_system}

\subsection{Equations of motion}
The mission design problem will be studied only in the cruise phase of the flight under the assumption that, if the flights have the same origin or the same destination, formation flight would be forbidden for security reasons during departure, climb, descent, and approach phases of flight. A simplified two-degree-of-freedom point variable-mass dynamic model is considered assuming that all aircraft are in the cruise phase. The motion is restricted to the horizontal plane at  cruise altitude over a spherical Earth model. A symmetric flight without sideslip has been considered and all the aircraft forces are supposed to be in the plane of symmetry of the aircraft. Wind effects have \textcolor{black}{also} been considered. All aspects associated to the rotational dynamics are neglected.  
The set of kinematics and dynamics Differential-Algebraic Equations (\textsf{DAE})  \cite{hull:2007:foafm} that \textcolor{black}{describe} each aircraft motion is
\begin{eqnarray} \label{eq:kinematics_dinamics_eqs}
\nonumber  \dot{\phi}(t) &   =  &   \dfrac{V(t) \cdot \cos \chi(t) + {V_{W_N}}(t)}{R_E + h},\\  
  	\dot{\lambda}(t) &   = &   \dfrac{V(t) \cdot \sin \chi(t) + {V_{W_E}}(t)}{ \cos \phi(t) \cdot \left( R_E + h \right) }, \\  
\nonumber  \dot{\chi} (t) &   =  &  \dfrac{L(t) \cdot \sin \mu(t)}{V(t) \cdot m(t)},\\ 
\nonumber   \dot{V}(t) &   =  &   \dfrac{T(t) - D(t)  }{m(t)}, 
\end{eqnarray}
where the state vector has five components, the two dimensional position variables, latitude and longitude, denoted by $\phi$ and $\lambda$, respectively, the heading angle $\chi$, the true airspeed $V$, and the mass of the aircraft $m$. In this set of equations, the control vector has three components, the thrust force $T$, the lift coefficient $C_L$, and the bank angle $\mu$. The normalized version of the \textsf{DAE} system (\ref{eq:kinematics_dinamics_eqs}) has been used in this paper.

The lift force is $L = qS C_L $, 
where $q = \frac{1}{2} \rho V^2$ is the dynamic pressure, $\rho$ is the air density, and $S$ is the reference wing surface area. The aerodynamic drag force is $D = qS C_D $, 
where $C_D$ is the drag coefficient. A parabolic drag polar $C_D = {C_D}_0 + K C_L^2$ has been assumed, where ${C_D}_0$ is the zero-lift drag component and $K$ is the induced drag coefficient. 
Other parameters and variables include the Earth radius $R_E$, the gravitational acceleration $g$, and the cruise altitude $h$. ${V_{W_E}}$ and ${V_{W_N}}$ are the components of the wind velocity vector in eastward and northward directions, respectively.
The mass flow rate equation to be added to the above set of \textsf{DAE} is
\begin{equation}
\dot{m}(t) = - T(t) \cdot  \eta(t), 
\label{mass_differential_eq} 
\end{equation}
where $\eta$ is the thrust specific fuel consumption. The Eurocontrol\textquotesingle s Base of Aircraft Data (\textsf{BADA}), version 3.6 \cite{eurocontrol:2013:umftboad}, has been used to determine it. In particular for jet engines the following equation applies
\begin{equation*}
\eta(t) = C_{f_1}\left( 1 + \frac{V(t)}{C_{f_2}} \right) ,
\label{eq:jet_engines_thrust}
\end{equation*}
where $C_{f_1}$ and $C_{f_2}$ are empirical thrust specific fuel consumption coefficients.

Thus, for aircraft $p$, the state vector will be  $x_p(t) = \left( \phi_p(t), \lambda_p(t), \chi_p(t), V_p(t), \right.$
$\left. m_p(t)\right) , \forall p \in  \left\lbrace 1, \ldots, N_a\right\rbrace $, and the control vector will be  $u_p(t) = \left( T_p(t), C_{L_p}(t), \right.$
$\left. \mu_p(t)\right) ,\forall p \in \left\lbrace 1, \ldots, N_a\right\rbrace $, with $N_a$ the number of aircraft involved in the mission design problem.

\subsection{Flight envelope}
Flight envelope constraints model aircraft performance limitations. These constraints involve, among others, flight altitude, load factor and airspeed, and can be stated in the form $\eta_l(t) \leq \eta(t) \leq \eta_u(t)$, where $\eta(t)$ is the state or control variable and $\eta_l(t)$ and $\eta_u(t)$ are the  minimum and maximum allowed values of the variable, respectively.
\textsf{BADA} 3.6 has been used to impose the following constraints
\begin{eqnarray} \label{flight_envelope_eqs}
\nonumber  V_{min}(t) \leq &  V_{CAS}(t) &  \leq   V_{MO},\\       
\nonumber  m_{min}  \leq   &  m(t)            &  \leq m_{max},\\       
\nonumber  			 &  M(t)            &  \leq M_{MO},\\        
            T_{min}(t) \leq         &  T(t)             &  \leq   T_{max}(t),\\ 
\nonumber   {C_L}_{min} \leq &  C_L (t) &  \leq {C_L}_{max},\\   
\nonumber  	-\mu_{max}  \leq&  \mu(t)&  \leq \mu_{max},      
\end{eqnarray}
where $V_{CAS}$ is the calibrated airspeed, which can be calculated as a function of the true airspeed \cite{eurocontrol:2013:umftboad}. $V_{min}$ and $V_{MO}$ are the minimum and maximum operating calibrated speeds, respectively.
$m_{min}$ and $m_{max}$ are the minimum and maximum aircraft masses, respectively. $M$ is the Mach number and $M_{MO}$ is the maximum operating Mach number.
$T_{min}$ and $T_{max}$ are the minimum and maximum available engine thrusts, respectively. ${C_L}_{min}$ and ${C_L}_{max}$ are the minimum and maximum lift coefficients, respectively. 
Finally, $\mu_{max}$ is the maximum bank angle set by the air navigation regulations in civil flight. 
Furthermore, the BADA aircraft model includes the equation $V_{min}(t) = C_{V_{min}} V_{s}(t)$, which relates $C_{V_{min}} $, the minimum speed coefficient, and $V_{s}$, the stall speed for each flight phase.

For turbofan propelled aircraft, assuming standard atmosphere conditions, the maximum thrust, $T_{max}$, is defined by the following empirical expression
\begin{equation*}
 T_{max}(t) =  C_{T_{cr}} \cdot C_{T_{C,1}} \cdot \left[  1 - \frac{H_p(t)}{C_{T_{C,2}}} + C_{T_{C,3}} \cdot H_p^2(t) \right],
 \label{eq:maximum_thrust}
\end{equation*}
where $C_{T_{cr}}$ is the maximum cruise thrust coefficient.
$C_{T_{C,1}}$, $C_{T_{C,2}}$, and $C_{T_{C,3}}$ are the empirical thrust coefficients and $H_p$ is the geopotential pressure altitude.

\subsection{Wind model}
\label{subsect:wind_model}
ERA-Interim \cite{Dee:2011ex}
is a third-generation global atmospheric reanalysis
produced by the European Centre for Medium-Range Weather Forecasts (\textsf{ECMWF}).
Climate reanalysis  consists in systematic  and coherent assimilation   of the global meteorological  data  constrained by  the available observations obtained from different sources during the period of reanalysis, such as radiosondes or aircraft,  and prior state estimates
using an invariant consistent assimilation scheme and an integrated meteorological model. The third generation reanalysis data refers to the last generation, which was developed in 2010s, and significantly enhanced the spatio-temporal resolution of the former ones.

Gridded data products, which can be found in the ERA-Interim database, include a wide variety of atmospheric parameters relevant to the mission design problem, such as air temperature, pressure, and wind at different altitudes and sea-surface temperature. In this study, the components of the wind velocity in eastwards and northwards directions have been used.

ERA-Interim is presented as a regular latitude-longitude grid, with a spatial resolution from $0.25^\circ \times 0.25^\circ$ to $3^\circ \times 3^\circ$  for the deterministic reanalysis, at 37 atmospheric levels corresponding to different pressure levels from 1000 to 1 hPa. The temporal coverage is four times daily, getting six-hour data.
In this study, the chosen pressure level has been 200 hPa,  corresponding to the cruise altitude.
The selected spatial resolution has been $0.5^\circ \times 0.5^\circ$. Finally, the selected day has been {April 30, 2019, at 12:00}.
To include atmospheric data in the hybrid optimal control problem used to solve the mission design problem, an analytic function that approximates the data must be determined. For this purpose, Radial Basis Functions (\textsf{RBF}) \cite{buhmann2003radial} have been used.

Given a set of known input data, a function approximation at any point is obtained constructing a linear space which depends on the relative position of the evaluated point with  respect to the known observations according to an arbitrary distance measure. 
\textsf{RBF} consider interpolating functions of the form
\begin{equation*}
F(x) = c_0 + \sum_{i=0}^{N-1} {c_c} \varphi(\parallel x - R_i \parallel), 
\label{eq:RBF}
\end{equation*}
where $x \in \mathbb{R}^N$ is the input data vector and $N$ is its dimension, $F(x)$ is the function approximation, $\varphi$ is the basis function, and $R_i $ is a vector containing the centers, which are the reference points of the basis functions.
$c_0$ and $c_c$ are constant coefficients calculated by \textsf{RBF} method, in particular, $c_0$ is the bias term.
 The Euclidean norm is considered to compute the distance between the input points and the centers. 
There are different options for the choice of the basis functions.  In this paper, Gaussian basis functions have been selected.

\begin{figure}[ht!]
\centering
\renewcommand{\figurename}{Fig.}
\includegraphics[width = 20em]{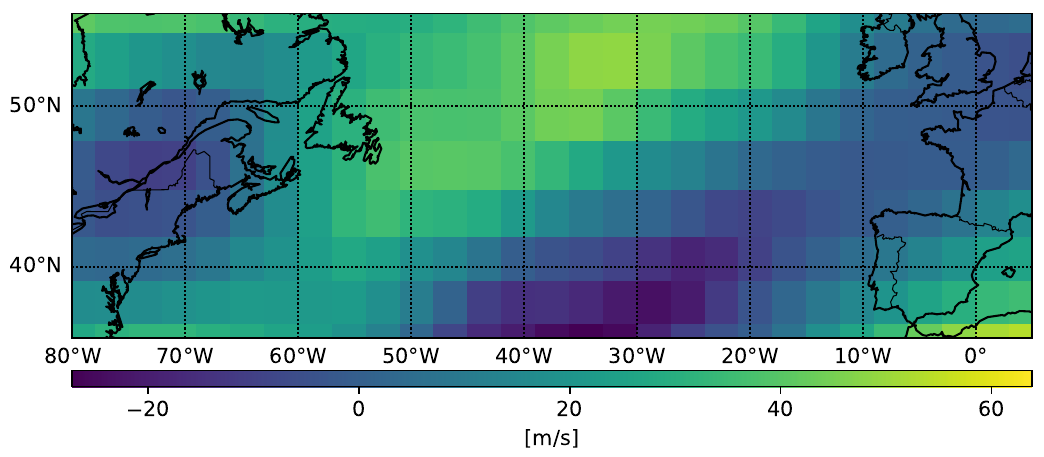}
\vspace{-1em}
\caption{Eastward wind speed on 30 April 2019, 12:00.}
\label{fig:WA_Ucomponent}
\end{figure}

\begin{figure}[ht!]
\centering
\renewcommand{\figurename}{Fig.}
\includegraphics[width = 20em]{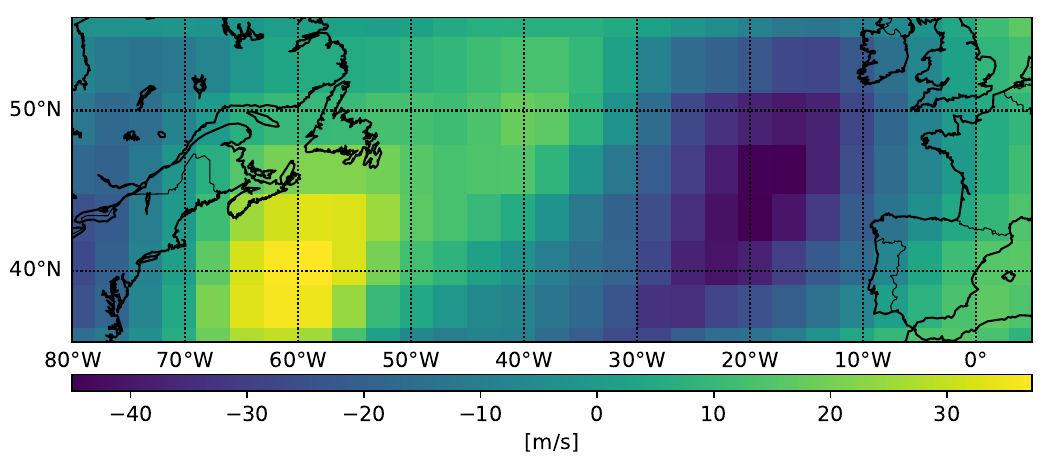}
\vspace{-1em}
\caption{Northward wind speed on 30 April 2019, 12:00.}
\label{fig:WA_Vcomponent}
\end{figure}

Eastwards and northwards components of the wind speed taken from ERA-Interim and then approximated by the \textsf{RBF} method are represented in Fig.~\ref{fig:WA_Ucomponent} and Fig.~\ref{fig:WA_Vcomponent}, respectively.

\subsection{Formation flight model}
\label{Formation_flight_model}
The \textsf{DAE} system consisting of Eq.~(\ref{eq:kinematics_dinamics_eqs}) and Eq.~(\ref{mass_differential_eq}) is applicable to any aircraft in the cruise phase of the flight. In solo flights, depending on aircraft type, the aerodynamic coefficients, the wingspan, and the reference surface, among others, will change. 
In formation flights, additionally, it is necessary to include a modification in Eq.~(\ref{mass_differential_eq}) for the trailing aircraft.

As previously mentioned, in extended formations, aircraft fly with a longitudinal separation of more than 10 wingspans but the potential fuel burn reduction for the trailing aircraft has a significant decrease beyond 20 wingspans separation. 
Therefore, in this paper, the formation benefits are considered negligible for separations of more than 20 wingspans and the constraint introduced to model the distance which leads to fuel burn reductions in the trailing aircraft, and in the intermediate one in case of three-aircraft formation, during formation flight is the following
\begin{equation}
10 b \leq \mathcal{D}_{pq}(t) \leq 20 b \label{eq_dist_ff},  \quad \forall p,q \in \lbrace1, \dots, N_a\rbrace, p < q, 
\end{equation}
where $b$ is the wingspan of the leader aircraft and $\mathcal{D}_{pq}(t)$ is the orthodromic distance between aircraft $p$ and $q$, at time $t$, given by the Haversine formula.

\begin{figure}[htb]
\centering
\renewcommand{\figurename}{Fig.}
\includegraphics[scale=0.4]{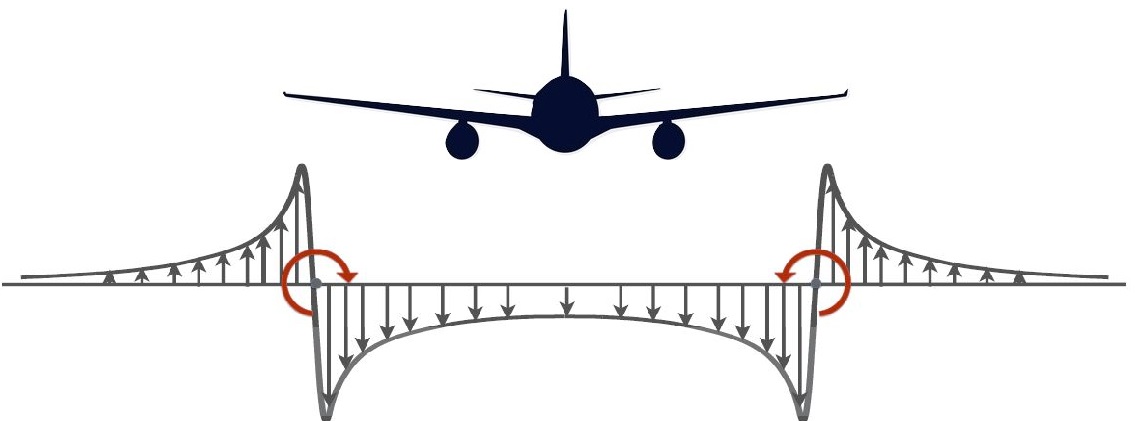}
\caption{Vertical component of the induced air velocity field (not to scale). Figure adapted from \cite{ning:2011:adrteff}.}
\label{fig:induced_downwash_upwash}
\end{figure}

Induced drag reduction in formation flight is due to the fact that any aircraft moving through a fluid generates lift by imparting a swirling motion to the air, which creates downwash and upwash fields. 
An aircraft flying in the upwash field created by another one will increase its apparent angle of attack,  getting an effective rotation of the lift vector resulting in a lower lift-dependent component of drag. In this way, a significant reduction in induced drag is obtained \cite{ningetal:2011:apoeff}. In Fig.~\ref{fig:induced_downwash_upwash} a schematic representation of the vertical component of the induced air velocity field is given.

A simplification has been introduced assuming that the reduction in the induced drag of the trailing and intermediate aircraft can be modeled as a percentage of the reduction in \textcolor{black}{fuel consumption}, $\mathcal{R}_\text{fuel}$. In this manner, the fuel saving factor can be directly incorporated in Eq.~(\ref{mass_differential_eq}) to obtain the mass flow rate equation for the trailing aircraft during formation flight
\begin{equation}  \label{mass_differential_eq_with_reduction}
\dot{m}(t) = - [1 - \mathcal{R}_\text{fuel}] \cdot T(t) \cdot  \eta(t) .
\end{equation}
It is quite difficult to estimate the drag reduction and, hence, the fuel burn reduction achieved during a formation flight, since there are many parameters to take into account. The number of aircraft involved in the formation, their type, size, and weight, the ability of the trailing aircraft to place itself in the optimum location of the wake, and the stresses and vibrations induced, among others, strongly influence the benefits of the formation. In this paper,  different fuel savings schemes for the trailing and the intermediate aircraft of a two- or three-aircraft formation have been considered. The adopted percentages are between  6\%  to 14\% of fuel savings for the trailing and intermediate aircraft, with respect to fuel consumption in solo fights.

\subsection{Objective {functional}}

The objective {functional} is a key factor in all trajectory optimization problems, but it is crucial in formation mission design problems, in which opposing optimality criteria, such as, fuel consumption and total mission time, must be considered.
On the one hand, formation flight can reduce the total fuel consumption and, on the other hand, in general, it requires increasing the total flight time to create the formation itself.
Therefore, formation mission design problems require to leave traditional fuel minimization models behind. 
In this paper, fuel consumption and flight time are combined in the objective functional using weighting coefficients.
The weighting coefficients have been chosen based on \cite{cook2004evaluating}, in which a detailed study on the costs incurred by airlines due to delays has been made.
 No other costs have been considered. 

The objective {functional} $J$ for $N_a$ aircraft, is defined as follows
\begin{equation}
\min \quad  J = \alpha_t  \sum_{p=1}^{N_a}  t_{{flight}_p}  +  \alpha_f \sum_{p=1}^{N_a}  {{{m}}_{f_p}}, 
\label{eq:obj_function}
\end{equation}  
where $\alpha_t$ and $\alpha_f$ are the time and the fuel burn weighting parameters, respectively, $t_{{flight}_p},  \forall p  \in \left\lbrace 1, \ldots, N_a\right\rbrace $, is the total flight time for aircraft $p$, 
and ${{{m}}_{f_p}} , \forall p  \in \left\lbrace 1, \ldots, N_a\right\rbrace $,  is the fuel consumption for aircraft $p$.
It is worth mentioning that the initial fuel loaded in an aircraft that is going to benefit from formation flight could be considerably lower than the amount of fuel loaded in the same aircraft but without reasonable chances to benefit from a formation flight. Nevertheless, in this paper, a more conservative approach has been adopted and all aircraft are loaded with the same amount of fuel as in solo flight, even if this fact will lead to lower benefits in terms of \textsf{DOC}.

As said before, in the approach to the formation mission design problem presented in this paper, the relative position of each aircraft is set in advance. Additionally, in three-aircraft formations, the only type of formation allowed is the in-line formation. A maximum detour time from the solo flight time has been set in order to 
perform a more realistic problem formulation according to the usual operational requirements.

\section{The Switched Optimal Control Problem}
\label{sect:the_switched_optimal_control_problem}
As said in the introduction, in this paper, a novel approach to the formation mission design problem is presented, which is formulated as an Optimal Control Problem (\textsf{OCP}) for a switched dynamical system.

In this section, following \cite{bengeaetal:2011:ocossveicocp}, the formulation of a Switched Optimal Control Problem (\textsf{SOCP}) is introduced for two-switched dynamical systems. The dynamical model of a two-switched dynamical system can be represented by
\begin{eqnarray} \label{eq:31.1}
\nonumber   \dot{x}_S(t) & =  & f_{v_{S}(t)}(t, x_S(t), u_S(t)), \\  
  x_S(t_I)  & =  & x_I \in \mathbb{R}^n, \\  
\nonumber     x_S(t_F)  & =  & x_F \in \mathbb{R}^n, \\  
\nonumber  v_S (t )  & \in &  \{0, 1\}, \;  t_I \leq t \leq t_F. 
\end{eqnarray}

\textcolor{black}{
For the sake of clarity, in the formulation of the \textsf{SOCP}, the sub-index ``$S$" has been added to the symbols that denote the state and control variables. }
In this model, the continuously differentiable vector fields, 
$f_0, f_1: \mathbb{R} \times \mathbb{R}^n \times \mathbb{R}^m \rightarrow \mathbb{R}^n$, 
specify the dynamics of the two possible modes of the system.
The control input $u_S(t) \in \Omega \subset \mathbb{R}^m$ is constrained 
to belong, at each time instant,  to the bounded and convex set 
$\Omega$. The binary variable $v_S (t ) \ignore{ \in \{0, 1\} }$ is the 
mode selection variable that identifies which of the two possible 
system modes, $f_0$ or $f_1$, is active. Thus, both of them, $u_S(t)$ and $v_S (t)$, 
can be regarded as control variables. 
The initial time $t_I$, the initial state $x_S(t_I)$, the final time $t_F$, 
and final state $x_S (t_F)$ are assumed to be
restricted to a boundary set $B$ as follows: $(t_I, x_S (t_I), t_F , x_S (t_F )) \in B = T_I \times B_I \times T_F  \times B_F  \subset \mathbb{R}^{2n+2}$.
The performance index of the \textsf{SOCP} has the following form
\begin{eqnarray}
&& J_S \left(t, x_S(t) , u_S(t) , v_S(t) \right) = g (t_F, x_F) + \int_{t_I}^{t_F}  F_{v_S(t)} (t, x_S (t), u_S (t)) dt,
\label{eq:31.2}
\end{eqnarray}
where the {endpoint cost} function $g$ is defined on a neighborhood of $B$, and $F_0$
and $F_1$ are real-valued continuously differentiable functions that
represent the {running} cost of operation of the system in each mode.

The \textsf{SOCP} is thus stated as follows
\begin{equation}
\min_{u_S \in \Omega, v_S \in \{0,1 \} } J_S \left(t, x_S(t) , u_S(t) , v_S(t) \right),
\label{eq_SOCP_definition}
\end{equation}
subject to Eq.~(\ref{eq:31.1}) and endpoint constraints $(t_I, x_S (t_I ), t_F , x_S (t_F)) \in B.$

\textcolor{black}{
In the formation mission design problem studied in this paper, not all the dynamic equations are subject to switches. More specifically, the equations of motion associated to the state variables ${\phi}$, ${\lambda}$, ${\chi}$, and ${V}$, do not switch whereas the aircraft\textquotesingle s mass flow rate equation
does when the aircraft joins or leaves a formation as a trailing or an intermediate aircraft. 
}

\textcolor{black}{
When an aircraft is flying solo or as a leader in the formation, it will obtain no benefits from the formation flight in terms of fuel savings and Eq.~(\ref{mass_differential_eq}) describes its mass flow. On the contrary, when an aircraft is flying as intermediate or trailing in the formation, it will obtain benefits from the formation flight and its mass flow will be described by Eq.~(\ref{mass_differential_eq_with_reduction}). Thus, $f_0$ represents Equations (\ref{eq:kinematics_dinamics_eqs}) and (\ref{mass_differential_eq}), whereas $f_1$ represents Equations (\ref{eq:kinematics_dinamics_eqs}) and (\ref{mass_differential_eq_with_reduction}). 
The functions $F_0$ and $F_1$, which appear in the running cost of Eq.~(\ref{eq:31.2}), do not exist in the formulation of this particular problem. 
}

\section{The Embedding Approach}
\label{sect:the_embedding_approach}
In this section, the embedding approach employed to transform the resulting \textsf{SOCP} into \textcolor{black}{an} Embedded Optimal Control Problem (\textsf{EOCP}) is presented. In this way, a smooth \textsf{OCP} without binary variables is obtained, reducing the computational complexity of finding the solution.
 
Following again  \cite{bengeaetal:2011:ocossveicocp}, Eq.~(\ref{eq:31.1}) can be rewritten as follows
\begin{eqnarray}
\dot{x}_E(t) &&= [1-v_E(t)] f_0(t, x_E(t), u_{E_0} (t))+ v_E(t) f_1 (t, x_E(t), u_{E_1} (t)).
\label{eq:31.3}
\end{eqnarray}

\textcolor{black}{
For the sake of clarity, in the formulation of the \textsf{EOCP}, the sub-index ``$E$" has been added to the symbols that denote the state and control variables. However, the components of each vector have not changed with respect to the \textsf{SOCP}. 
}

Equations (\ref{eq:31.1})  and (\ref{eq:31.3}) formally coincide under the conditions $v_{E}(t) = v_S(t)$ and $u_{E_0}(t) = u_{E_1}(t) = u_S(t)$.
However, in Eq.~(\ref{eq:31.3}),  $v_{E}(t) \in [0,1]$, while in Eq.~(\ref{eq:31.1}), $v_S(t) \in \{0,1\}$, 
which is the key feature of the embedding  approach \cite{bengeaanddecarlo:2005:ocoss, bengeaetal:2011:ocossveicocp}.
This method is based on solving the \textsf{SOCP} defined above, using only continuous control variables  $v_{E}(t)  \in [0,1] $, 
$u_{E_0}(t), u_{E_1}(t) \in \Omega$, in which the performance index  of the \textsf{SOCP} {from Eq.}(\ref{eq:31.2}) is rewritten as 
\begin{eqnarray}
&&J_E \left(t, x_E(t), u_{E_0}(t), u_{E_1}(t), v_E(t) \right) = g ( t_F, x_E(t_F))  + 
 \nonumber \\
&&+ \int_{t_I}^{t_F} \left(  [1- v_E (t)] F_0 (t, x_E(t), u_{E_0} (t)) + v_E(t)  F_1(t, x_E(t), u_{E_1} (t) ) \right) dt. 
\label{eq:31.4}
\end{eqnarray}

The \textsf{EOCP} is thus formulated as
\begin{equation}
\min_{u_{E_0}, u_{E_1} \in \Omega, v_E \in [0,1] } J_E \left(t, x_E(t), u_{E_0}(t), u_{E_1}(t), v_E(t) \right),
\label{eq_EOCP_definition}
\end{equation}
subject to Eq.~(\ref{eq:31.3}) and endpoint contraints $(t_I, x_E (t_I ), t_F , x_E (t_F)) \in B$, which is an \textsf{OCP} without binary variables. Therefore, standard numerical optimal control techniques can be applied to solve
it.

It has been shown in \cite{bengeaanddecarlo:2005:ocoss} and \cite{bengeaetal:2011:ocossveicocp} that, 
once a solution of the \textsf{EOCP}  has been obtained,
either the solution is of the switched type, that is, 
$v_{E}$ takes only the values $0$ and $1$, or suboptimal trajectories of
the \textsf{SOCP}  can be constructed that can approach the value of the cost
for the \textsf{EOCP}  arbitrarily closely, and satisfy the boundary conditions
within $\epsilon$, with arbitrary $\epsilon > 0$. A thorough discussion
about the relationship between the solutions of the \textsf{SOCP} and the \textsf{EOCP} can be
found in \cite{bengeaanddecarlo:2005:ocoss} and \cite{bengeaetal:2011:ocossveicocp}. 

\subsection{Specification of the Embedded Optimal Control Problem}
\label{subsect:specification_of_the_embedded_optimal_control_problem} 
\textcolor{black}{
As said in Section III.A, only the aircraft\textquotesingle s mass flow rate equation is subject to switches.}
To model this switching, equations (\ref{mass_differential_eq}) and (\ref{mass_differential_eq_with_reduction}) are combined using the selection variable $v_{E_p}(t), \; p \in \left\lbrace 1, \ldots, N_a\right\rbrace $ as follows
\begin{eqnarray}  \label{mass_differential_eq_var_decision}
 \nonumber & \dot{\phi}_p(t) &  = \quad  \dfrac{V_p(t) \cdot \cos \chi_p(t) + {V_{W_{N}}}(t)}{R_E + h_p},\\  
\nonumber   &\dot{\lambda}_p(t) & =  \quad \dfrac{V_p(t) \cdot \sin \chi_p(t) + {V_{W_{E}}}(t)}{ \cos \phi_p(t) \cdot \left( R_E + h_p \right) }, \\  
&\dot{\chi}_p (t) & =  \quad  \dfrac{L_p(t) \cdot \sin \mu_p(t)}{V_p(t) \cdot m_p(t)},\\ 
\nonumber   & \dot{V}_p(t) & =   \quad \dfrac{T_p(t) - D_p(t)  }{m_p(t)}, \\
\nonumber  & \dot{m}_p(t) & =   \quad   (1 - v_{E_p}(t)) \cdot [- T_p (t) \cdot  \eta_p (t)] +  \\
\nonumber  & &   \hspace{0.75cm}+ \;v_{E_p}(t) \cdot [- (1 -\mathcal{R}_{\text{fuel}_p}) \cdot T_p(t) \cdot  \eta_p(t)],
\end{eqnarray}
\textcolor{black}{
$\forall p \in   \left\lbrace 1, \ldots, N_a\right\rbrace$, where $v_{E_p}(t) = 1$ corresponds to Formation Flight in \textcolor{black}{an} in-line formation (\textsf{FF}) and $v_{E_p}(t) = 0$ to Solo Flight (\textsf{SF}).}

Therefore, each aircraft has different flight modes, namely \textsf{SF}
and \textsf{FF} in different positions inside a formation. Their combination is represented by the discrete state of the switched system which models their joint dynamic behaviour.

However, in the mission design problems, the switching control variables $v_{E_p}(t), \; p \in \left\lbrace 1, \ldots, N_a\right\rbrace $  depend on the state variables, in particular on the geographical coordinates of the aircraft. Indeed, not all transitions between discrete states are possible at any time, since while the distance between aircraft is larger than 20 wingspans, \textsf{FF} can not take place. Furthermore, the switching control variables depend on one another. For instance, in a two-aircraft formation, if one aircraft is selected as the trailing aircraft the other one is forced to be the leader and vice versa. Similarly, in a three-aircraft formation, if one position in the formation is assigned to one aircraft, the other ones must fly in another position.

Thus, the switched system that represents the joint dynamic behavior of all the aircraft involved in the mission design problem should follow some rules that are specified using logical constraints usually expressed in disjunctive form. The method employed to model logical constraints in disjunctive form in the \textsf{OCP} used to solve the mission design problem will be described in the next section.

\section{Logical Constraints Modeling}
\label{sect:treatment_of_the_logical_constraints} 
To transform logical constraints in disjunctive form into inequality
and equality constraints that involve only continuous auxiliary variables, the approach proposed in \cite{weietal:2008:ocorswlcatupp} has been employed.

It has been shown in \cite{cavalieretal:1990:maiptatpc} that every
Boolean expression can be transformed into Conjunctive Normal Form (\textsf{CNF}). Hence, there is no loss of generality in considering  that any logical constraint can be rewritten as a \textsf{CNF} expression
\begin{equation}
Q_1 \wedge Q_2 \wedge \ldots \wedge Q_n,
\label{logical_formula_conjunction} 
\end{equation}
in conjunctive form, where the symbol $\wedge$ denotes the $``$and$"$ operator and 
\begin{equation}
Q_i = P_i^1 \vee P_i^2 \vee \ldots \vee P_i^{m_i}, \; \forall i \in \{ 1,2,\ldots, n \},
\label{logical_formula_disjunction}
\end{equation}
where $Q_i$ is expressed in disjunctive form, in which the symbol $\vee$ denotes the $``$or$"$ operator and the  proposition $P_i^{j}$ is either $X_i^{j}$ or $\neg X_i^{j}$. The term $X_i^{j}$  can be either True or False and the symbol $\neg$ represents the $``$not$"$ operator.

Term $X_i^{j}$ represents statements such as $``$$ \mathcal{D}_{12}(t) \geq 20b \;$$"$. Therefore, $P_i^{j}$ takes the form 
\begin{equation}
P_i^{j} \equiv \left\lbrace   g_i^j \left( x(t)\right)    \leq 0  \right\rbrace ,
\label{logical_formula_inequality}
\end{equation}
$\forall i \in \{ 1,2,\ldots, n \}$, $\forall j \in \{ 1,2,\ldots, m_i \}$, where $ g_i^j: \mathbb{R}^{n_x} \rightarrow \mathbb{R}$  is assumed to be a $\mathcal{C}^1$ function.

To incorporate logical constraints in a \textsf{OCP}, they must be converted into a set of equality or inequality constraints. Additionally, if binary variables are not used, the combinatorial complexity of integer programming is avoided.
Conjunctions from Eq.~(\ref{logical_formula_conjunction}) 
can be straightforwardly included into a \textsf{OCP} since they are equivalent to 
$
Q_i, \forall i \in \{ 1,2,\ldots, n \}$.
To be able to transform the disjunctions into a set of inequality constraints, a continuous variable $\alpha_i^j(t)\in [0,1]$ is defined and associated with each $P_i^j$ in Eq.~(\ref{logical_formula_inequality}). Thus, Eq.~(\ref{logical_formula_disjunction}) can be rewritten as
\begin{eqnarray}
&&   \alpha_i^j (t)  \cdot g_i^j \left( x(t) \right)  \leq 0, \nonumber  \\ 
\text{and} && 0 \leq \alpha_i^j (t) \leq 1, \label{logical_formula_disjunction_inequality_set} \\  
\text{and} && \sum_{j=1}^{m_i} \alpha_i^j (t)  =1, \nonumber
\end{eqnarray}
$\forall i \in \{ 1,2,\ldots, n \}, \forall j \in \{ 1,2,\ldots, m_i \}$. The first constraint in (\ref{logical_formula_disjunction_inequality_set}) shows that when $\alpha_i^j(t)=0$ the constraint $g_i^j \left( x(t)\right)   \leq 0$ is not forced to be fulfilled. In contrast, if $0 < \alpha_i^j (t) \leq 1$, then  $ \alpha_i^j(t)  \cdot g_i^j \left( x(t)\right)   \leq 0$, and, consequently, the constraint $g_i^j \left( x(t)\right)    \leq 0 $ is forced to be fulfilled. The last equation in (\ref{logical_formula_disjunction_inequality_set}) guarantees that at least one of the propositions $P_i^j$ holds.

In the next section, this approach will be applied to the formation mission design problem in which logical constraints involve the longitudinal distance between aircraft.

\subsection{Specification of the Logical Constraints}
\label{subsect:specification_of_the_logical_constraints_modeling} 
As stated in Eq.~(\ref{eq_dist_ff}), in this paper, it has been considered that the great-circle distance between aircraft in \textsf{FF} is between 10 and 20 wingspans, and is greater than 20 wingspans in case of \textsf{SF}.  Hence, the distance between two aircraft, $\forall p,q\in \lbrace 1, \dots, N_a \rbrace, \; p < q$, should satisfy
\begin{eqnarray}
&  \; \;10 b  \leq  & \mathcal{D}_{pq}(t) \; \; \leq 20 b, 
\label{eq_distancias_or} \\  
\text{or}  &                   & \mathcal{D}_{pq}(t)\; \; \geq 20 b. \nonumber
\end{eqnarray}
Using simple Boolean algebra,  Eq.(\ref{eq_distancias_or}) can be rewritten as
\begin{eqnarray}
&\; \;  10 b  \leq   \mathcal{D}_{pq}(t),\label{eq_distancias_or2} 
\\  
\text{and}   & \quad \mathcal{D}_{pq}(t)\; \; \leq 20 b   \quad 
\text{or}          \quad             \mathcal{D}_{pq}(t)\; \; \geq 20 b. \nonumber
\end{eqnarray}

To transform the disjunctions into a set of equality and inequality constraints, a new variable is defined. Thus, in two-aircraft mission design problems, only a continuous variable related to the longitudinal distance $\mathcal{D}_{12}$ should be defined, whereas, in three-aircraft mission design problems three continuous variables must be introduced, which are related to the longitudinal distances $ \mathcal{D}_{12},  \mathcal{D}_{13},  \mathcal{D}_{23}$. Thus, calling the new continuous variable $\alpha_{pq}(t)$ $\in [0,1]$, $\forall p,q \in \lbrace 1, \dots, N_a \rbrace, \; p < q $, the logical constraints in Eq.(\ref{eq_distancias_or2}) can be rewritten as
\begin{eqnarray}
&&  \alpha_{pq}(t) \left( \mathcal{D}_{pq}(t)- 20b \right)  \leq 0,  \nonumber  \\ 
  \text{and}  &&  \left(1- \alpha_{pq}(t) \right)  \left( 20b - \mathcal{D}_{pq}(t)\right)  \leq 0,   \nonumber \\
\text{and}  &&   0 \leq \alpha_{pq}(t) \leq 1, \label{eq_logical_constraints} \\
\text{and}  &&  \mathcal{D}_{pq}(t) \geq 10b. \nonumber  
\end{eqnarray}

For aircraft $p$ and $q$ $\in \lbrace 1, \dots, N_a \rbrace, \; p < q $, if $ \mathcal{D}_{pq}(t) > 20b$, $\alpha_{pq}(t)$ must be zero in order to satisfy the first constraint in (\ref{eq_logical_constraints}). On the contrary, when $\mathcal{D}_{pq}(t)  \leq 20b$, the variable $\alpha_{pq}(t)$ must be equal to one to fulfill the second constraint in the above set of equations. Formation flight benefits will be achieved only when $\alpha_{pq}(t) = 1$. In any case, the last constraint in (\ref{eq_logical_constraints}) ensures that the longitudinal distance between aircraft will always be greater than $10b$,  the minimum safety distance.

It can be observed that nature 
of the variables $v_{E_p}(t)$, $\forall p \in \lbrace 1, \dots, N_a \rbrace$, and $\alpha_{pq}(t)$, $\forall p,q \in \lbrace 1, \dots, N_a \rbrace, \; p < q$, is similar. While the variable $v_{E_p}(t)$ selects equations of the dynamic model for aircraft $p$, the variable  $\alpha_{pq}(t)$ selects constraints. Although, for the sake of clarity, a different notation has been used in the previous sections to describe them, the variables  $\alpha_{pq}(t)$ and  $v_{E_p}(t)$  are directly related by the following equivalence
\begin{equation}
\forall p,q \in \lbrace 1, \dots, N_a \rbrace, \; p < q :  \quad  \alpha_{pq}(t) = v_{E_p}(t), \label{eq:equivalence_variables_continuas}
\end{equation}
where aircraft $p$ is supposed to be the trailing aircraft. Therefore, the variable $v_{E_q}(t)$, which selects equations of the dynamic model for aircraft $q$, must be zero regardless of the value of variable $\alpha_{pq}(t)$, as aircraft $q$ will be the leader one and, hence, it will achieve no benefits from the formation.

It is important to notice that the equivalence in Eq.~(\ref{eq:equivalence_variables_continuas}) between the continuous variables obtained from the embedding approach, used to model the switched dynamical system and transform the logical constraints in disjunctive form into inequality and equality constraints, 
 need not necessarily to be fulfilled in others \textsf{EOCP}s.

\section{Solution of the Switched Optimal Control Problem}
\label{sect:solution_of_the_switched_optimal_control_problem}
Once the \textsf{SOCP} has been reformulated as an \textsf{EOCP}  and the logical constraints in disjunctive form have been transformed into inequality and equality constraints without binary variables using the embedding approach, the resulting optimal control problem has been transcribed into a NonLinear Programming (\textsf{NLP}) problem
using a PeudoSpectral Method (\textsf{PSM}). More specifically, the Radau \textsf{PSM} has been employed \cite{Garg2017AnOO}.

The set of constraints of the resulting \textsf{NLP} problem includes the system constraints that correspond to the differential constraint (\ref{mass_differential_eq_var_decision}) and the discretized versions of the other constraints of the  \textsf{OCP}, including  the state and control constraints
(\ref{flight_envelope_eqs}),  the algebraic constraints (\ref{eq_logical_constraints}), 
and the boundary conditions.
The resulting  \textsf{NLP} problem has been modeled using Pyomo \cite{hartetal:2012:pomip}, an open-source Python-based software package for modeling complex optimization problems,  and solved using the {Interior Point OPTimizer} (\textsf{IPOPT}) solver, an open-source software package suitable for large-scale nonlinear optimization. It implements an interior-point line-search filter method and can be used to solve general \textsf{NLP}  problems \cite{wachter2006implementation}.

To solve smooth \textsf{OCP}, \textsf{PSM} {are} more accurate than traditional
collocation methods \cite{betts2010practical}. However, as they are based on global polynomials with predetermined location of the nodes, 
their use in solving \textsf{SOCP} can give rise to numerical difficulties since the location of the switches may not coincide with the preallocated nodes.
Hence, to solve these difficulties, the concept of pseudospectral knots was introduced in \cite{Ross2004PseudospectralKM}. 
{The knots are some special nodes which enable the global time interval to be divided into several time subintervals to apply \textsf{PSM} over each subinterval, }
allowing the use of \textsf{PSM} to solve  \textsf{SOCP}, such as formation mission design problems. Therefore, a pseudospectral knotting method has been used to solve the \textsf{OCP} resulting from the \textsf{EOCP}.

\subsection{The Pseudospectral Method}
\label{sect:pseudospectral_method}

\textsf{PSM} are a type of 
collocation methods which use global polynomials to parameterize the state variables. The nodes are generally obtained from a Gaussian quadrature, which is used to collocate the differential-algebraic equations. The use of global polynomials, instead of using piecewise-continuous polynomials as interpolant between prescribed subintervals, makes \textsf{PSM} more efficient and simpler than other direct methods \cite{fahroo2000spectral}.

Due to orthogonality properties, the first step to introduce \textsf{PSM} is to transform the flight time interval of each flight $p$, defined as
$t_p \in I_p = [t_{I_p}, t_{F_p}], \; p \in \lbrace 1, \dots, N_a \rbrace$, into a new set of time variables, $\tau_p\in [-1,+1], \;  p \in \lbrace 1, \dots, N_a \rbrace$.
The transformation that relates the original time interval of each aircraft, $t_p$, 
is the following 
\begin{equation*}
t_p = \dfrac{\left( t_{F_p} - t_{I_p} \right)\tau_p   \;  +  \; \left( t_{F_p}+ t_{I_p}\right)}{2} .
\label{time_transformation}
\end{equation*}

In trajectory optimization problems, the three most commonly used sets of collocation points  are Legendre-Gauss (LG), Legendre-Gauss-Radau (LGR), and Legendre-Gauss-Lobatto (LGL) points. Notice that, although in all these methods 
the discretization points used in the \textsf{NLP} problem must lie in the $[-1,+1]$ interval, this is not required for the collocation points.
Indeed, depending on the set of collocation points chosen, both, one or none endpoints of the interval $[-1,+1]$ can be excluded in the approximation of the control variables.
In this paper, a \textsf{PSM} based on collocation at Flipped Radau points have been used \cite{Garg2017AnOO}. The Flipped Radau collocation points lie in the interval $(-1,+1]$.

Additionally, the \textsf{EOCP} stated to solve the formation mission design problem has some peculiarities which led to the use of the
knotting method in the \textsf{PSM} that will be outlined in the next section.

\subsection{The Pseudospectral Knotting Method}
\label{sect:knotting_method}
The obtained \textsf{EOCP} is a nonsmooth problem with switches in the state and control variables, and in the dynamic constraints. 
In this \textsf{EOCP}, the initial and final times are supposed to be different for each \textcolor{black}{flight} and, consequently, the nodes of each flight do not match. However, to enforce the separation constraints between aircraft, it is necessary to have the same nodes for all the flights involved, at least, during the formation flight. The pseudospectral knotting method solves this necessity in the following way. 
They generalize the spectral patching method by the exchange of information across the patches between two standard \textsf{PSM} in the form of switching conditions \cite{Ross2004PseudospectralKM}. These switching conditions are localized at the so-called knots of the problem. Hence, the knots represent double nodes, coinciding with the last node of one subinterval or patch and with the first node of the following subinterval. Different kind of knots can be defined, such as soft or hard knots, or free or fixed knots. 

Following \cite{Ross2004PseudospectralKM},  
in the formation mission design problems solved in this paper, only two interior soft knots are considered, independently of the number of aircraft involved in the formation.
Besides, initial-time and final-time conditions are considered as part of the general framework of knots. In particular, the knots corresponding with the initial-time and final-time are defined as hard knots because they are intrinsic to the problem formulation.

\begin{figure}[htb]
\centering
\renewcommand{\figurename}{Fig.}
\includegraphics[scale=0.5]{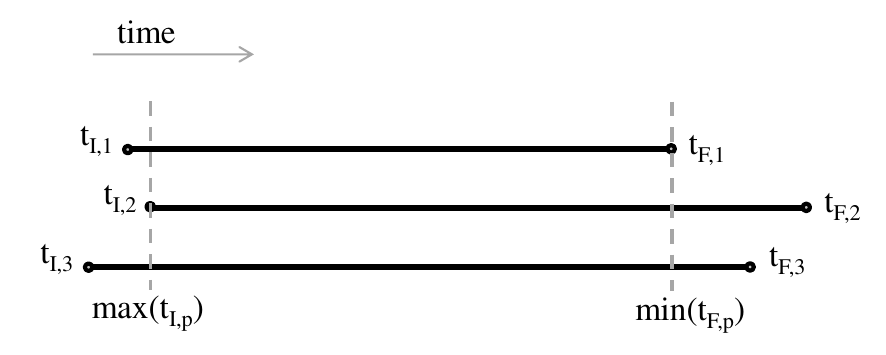}
\caption{Schematic flight times representation in a three-aircraft mission design problem.}
\label{fig:times_representation_3aircraft}
\end{figure}

The interior knots are defined as free knots but constrained by the condition that the time position of the first interior knot, $t_{k_1}$, should be greater than the initial times of all flights, $\text{max}(t_{I,p})$, and the time position of the second one, $t_{k_2}$, should be smaller than the final times of all flights, $\text{min}(t_{F,p})$, as schematically shown in Fig.~\ref{fig:times_representation_3aircraft}. 

On this basis, three different time intervals, $I^1$, $I^2$,  and $I^3$,  have been considered for each flight $p \in \lbrace 1, \dots, N_a \rbrace $:
\begin{itemize}
\item[•] \text{$I^1$: before the first interior knot},  $t_p \in [t_{I,p} ,\; t_{k_1})$, formation is not allowed.
\item[•] \text{$I^2$: between the two interior knots},  $t_p \in [ t_{k_1}, \; t_{k_2}]$, formation is allowed.
\item[•] \text{$I^3$: after the second interior knot}, $t_p \in ( t_{k_2}, \; t_{F,p}]$, formation is not allowed.
\end{itemize}

\begin{figure}[htb]
\centering
\renewcommand{\figurename}{Fig.}
\includegraphics[scale=0.3]{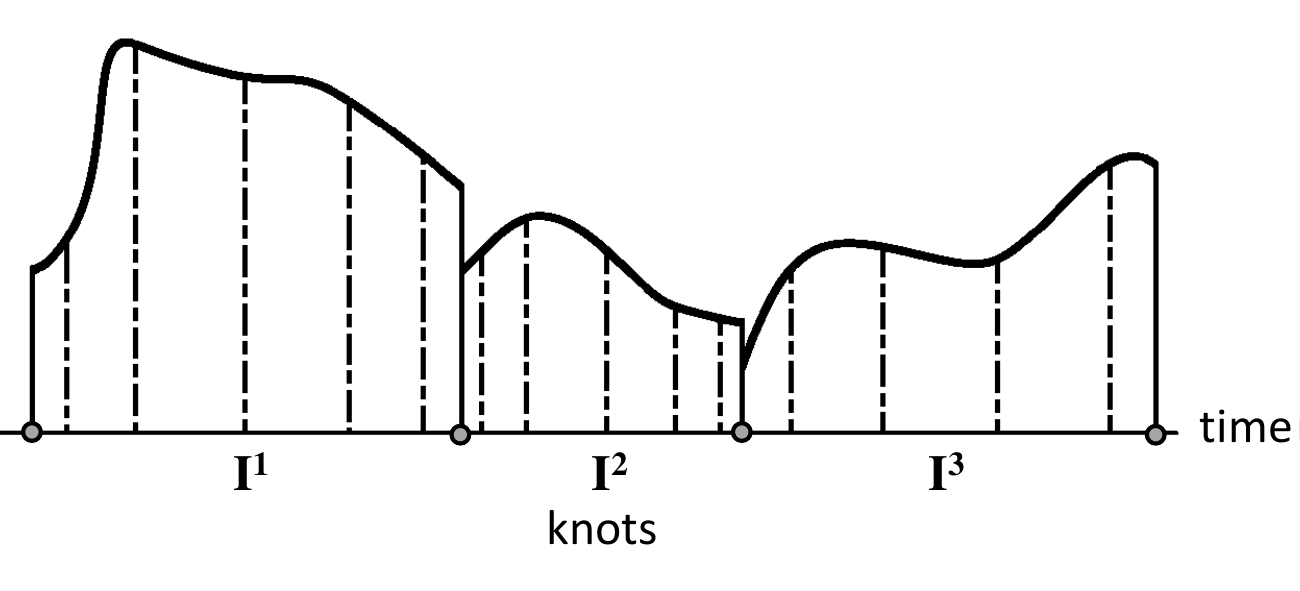}
\caption{Knots and nodes. Figure adapted from \cite{Ross2004PseudospectralKM}.}
\label{fig:knots_nodes}
\end{figure}

In Fig.~\ref{fig:knots_nodes} a schematic representation of the three different time intervals, $I^1$, $I^2$,  and $I^3$, and the knots and nodes of the problem are represented. The vertical lines correspond to the times of each node, and the big points represent the knots.

It is easy to see that, in two-aircraft formation design problems, the formation can 
only happen if both aircraft are flying, in other words, during the overlapped flight time \cite{xuetal:2014:arofff}.
Thus, the formation is only allowed during the second time interval, $I^2$, and the collocation points of all the aircraft are forced to be the same in this interval. Notice that, as mentioned before, the two interior knots of the problem should be greater than $\text{max}(t_{I,p})$ and lower than $\text{min}(t_{F,p})$, but are not fixed to be equal to those values. Hence,  the time value of the interior knots should also be determined in the solution. The same assumption has been taken for three-aircraft formations.

In addition to the numerical difficulties associated to the nonsmoothness of the mission design problem, the Gibbs phenomenon should also be mentioned.  This phenomenon appears when a nonsmooth function is approximated with some smooth functions  \cite{Ross2004PseudospectralKM}, and the pseudospectral knotting method helps to deal with it.

\section{Numerical Results}
\label{sect:numerical_results}

To show the effectiveness of the methodology for the formation mission design problem described in the previous sections, the following numerical experiments have been carried out: 

\begin{itemize}
\item Experiment A: two-aircraft transoceanic mission design with one flight\textquotesingle s departure time free.
\item Experiment B: two-aircraft transoceanic mission design with delays in the departure times.
\item Experiment C: three-aircraft transoceanic mission design with different fuel savings schemes.
\end{itemize}

All the experiments involve transoceanic eastbound flights. Wind data from the ERA-INTERIM reanalysis database \textcolor{black}{has} been used.
As mentioned above, in this paper only the cruise phase has been modeled, the rest of the flight phases are neglected. Thus, the initial and final cruise phase locations of the cruise phase of the flights have been assumed to be the latitudes and longitudes of the departure and arrival airports of each flight at  cruise altitude.  
Airbus A330-200 aircraft \textsf{BADA} models have been considered  for each flight. 
The initial mass of each aircraft is assigned as well as the initial and final velocities, which have been set at typical cruise values for the selected aircraft models. The initial heading angle is set to the initial heading angle of the orthodromic path between the initial and final locations of each flight.
Arrival times are left free and a constraint on the maximum temporal deviation from the scheduled arrival time of 45 minutes has been introduced for each flight.
The time and fuel burn weighting parameters, $\alpha_t$ and $\alpha_f$, in the objective functional (\ref{eq:obj_function}) have been set to 0.3 and 0.7, respectively \cite{cook2004evaluating}.

The numerical experiments have been conducted on a 3.6 GHz Intel Core i9 computer with 32 GB RAM. 
The computational times reported in this section include both the time required to generate the warm-start solution and the time to find the optimal solution of the problem.
The warm-start solution is a feasible solution of the problem generated by the NLP solver from initial guesses of the solution. The initial guesses for the latitude, the longitude and the heading angle have been generated using the orthodromic between departure and arrival locations of each flight. 
Typical cruise velocity and fuel consumption of the aircraft model selected have been used to generate the initial guesses of the velocity and the mass of each aircraft during the flight, respectively.

\subsection{Experiment A: two-aircraft transoceanic mission design with one flight\textquotesingle s departure time free}
Experiment A involves two transoceanic eastbound flights, Flight 1 and Flight 2, with given fuel savings for the trailing aircraft and boundary values of the state variables. Flight 1 and Flight 2 are operated by Aircraft 1 and Aircraft 2, respectively.
In the \textsf{SOCP} stated to solve the mission design problem, the departure time of Flight 1 is fixed, whereas the departure time of Flight 2 is left free. In case of formation, Aircraft 1 will be the trailing one and Aircraft 2 will be the leader.
{Solving this problem entails to decide which mode of flight, i.e. formation or solo flight, is optimal,  the optimal trajectories of the aircraft and the optimal departure time of Flight 2.}
Additionally, in case of formation flight, 
the rendezvous and splitting locations and times must also be determined. 
To check the optimality of the obtained solution, a  comparison between formation flight and solo flight results has been carried out.

The two transoceanic flights considered in this experiment have the following departure and arrival locations: 
\begin{itemize}
\item Flight 1: New York (\texttt{JFK}) - Madrid (\texttt{MAD}).
\item Flight 2: Montreal (\texttt{YUL}) - London (\texttt{LHR}).
\end{itemize}

The fuel burn savings for the trailing aircraft is set to 10\%.
The departure time of Flight 1 has been set to 10:15.  
The boundary conditions for the state variables of each aircraft are listed in Table \ref{table:boundary_conditions}.

\begin{table}[ht!]
\centering
\caption{Experiment A: boundary conditions.}
\medskip

\begin{tabular}{c c c c c}
\multicolumn{1}{ c }{\textbf{Symbol}} & \textbf{{Units}}& \textbf{{Flight 1}}             & \textbf{{Flight 2}}         \\ \hline
 $\phi_I$ & [deg] & 40.64  &  45.47 \\
 $\phi_F$ &  [deg] & 40.48  &  51.47  \\
 $\lambda_I$ &  [deg] & -73.78  & -73.74\\ 
 $\lambda_F$ &  [deg] & -3.57  &  -0.45  \\ 
 $\chi_I$ &  [deg] & 66.51    &  55.70   \\
 $V_I$ &  [m/s] & 240   &   240 \\
 $V_F$ &  [m/s] & 220   &   220 \\
 $m_I$ &  [kg] & $220 \, 000$    &  $215 \, 000$ \\ \hline
\end{tabular}
\label{table:boundary_conditions}
\end{table}

As mentioned above, the formation configuration is selected in advance. In case of formation flight, Aircraft 1, associated to Flight 1, is forced to be the trailing aircraft. Hence, Aircraft 2, associated to Flight 2, will be the leader one.

\begin{figure*}[htb]
\centering
\renewcommand{\figurename}{Fig.}
\includegraphics[scale=0.75]{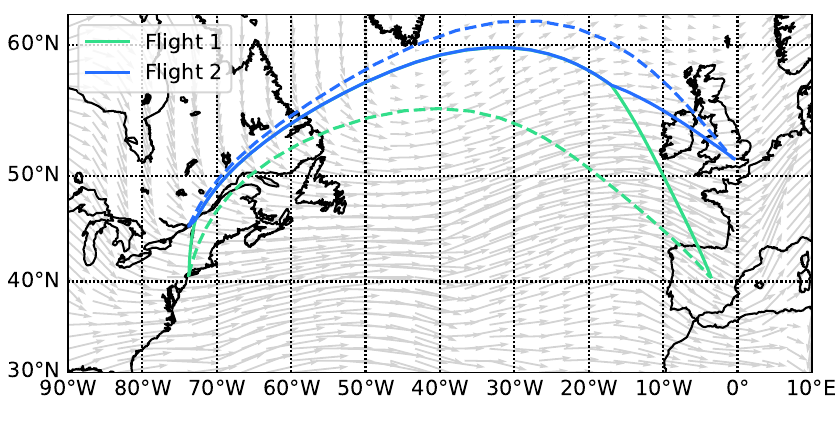}
\caption{Experiment A: solo and formation flight routes.}
\label{fig:baseline_wind_map}
\end{figure*}

In the optimal solution, formation flight has been selected as the optimal solution and the obtained optimal departure time for Flight 2 has been 10:50, 35 minutes after the departure of Flight 1. The obtained optimal formation flight routes are represented in Fig.~\ref{fig:baseline_wind_map},  
together with the solo flight routes and the wind field. 
In this figure, dashed lines are used to plot solo flight routes and 
solid lines are used to represent formation flight routes. Green and blue lines represent Flight 1 and Flight 2, respectively.

Flight 1 starts 35 minutes before Flight 2, and due to this difference in the departure times, the optimal solution implies a large initial detour of Flight 1 with respect to solo flight, in order to  join Flight 2 as soon as possible, thereby, the rendezvous point is very close to the departure location of Flight 2.

\begin{figure*}[ht!]
\centering
\renewcommand{\figurename}{Fig.}
\hspace{-2mm}
\subfigure[Mass]{\includegraphics[width=55mm]{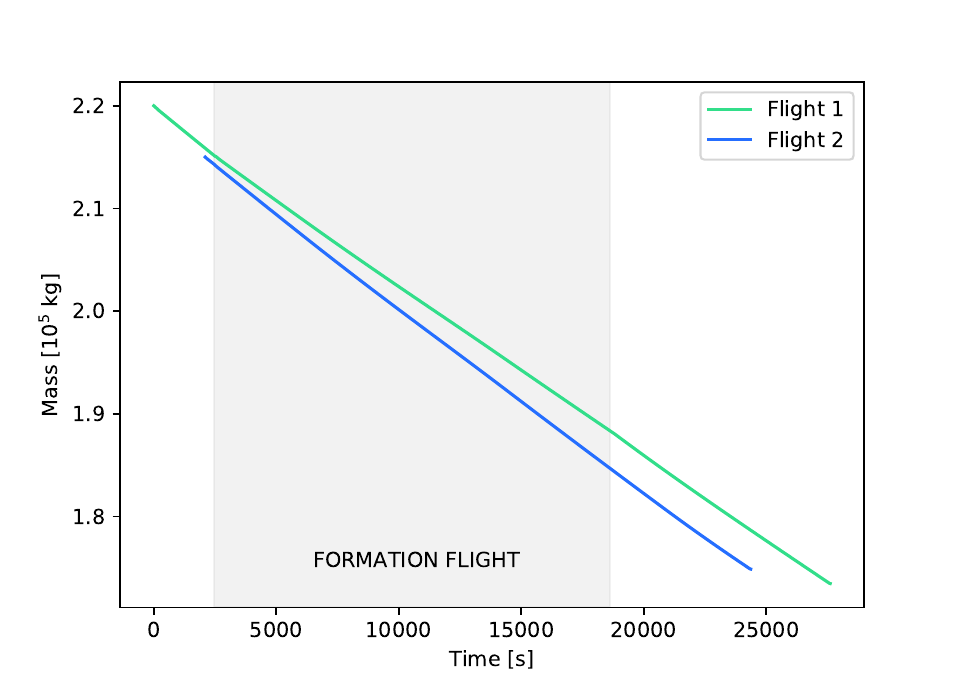}}
\hspace{-2mm}
\subfigure[Heading angle]{\includegraphics[width=55mm]{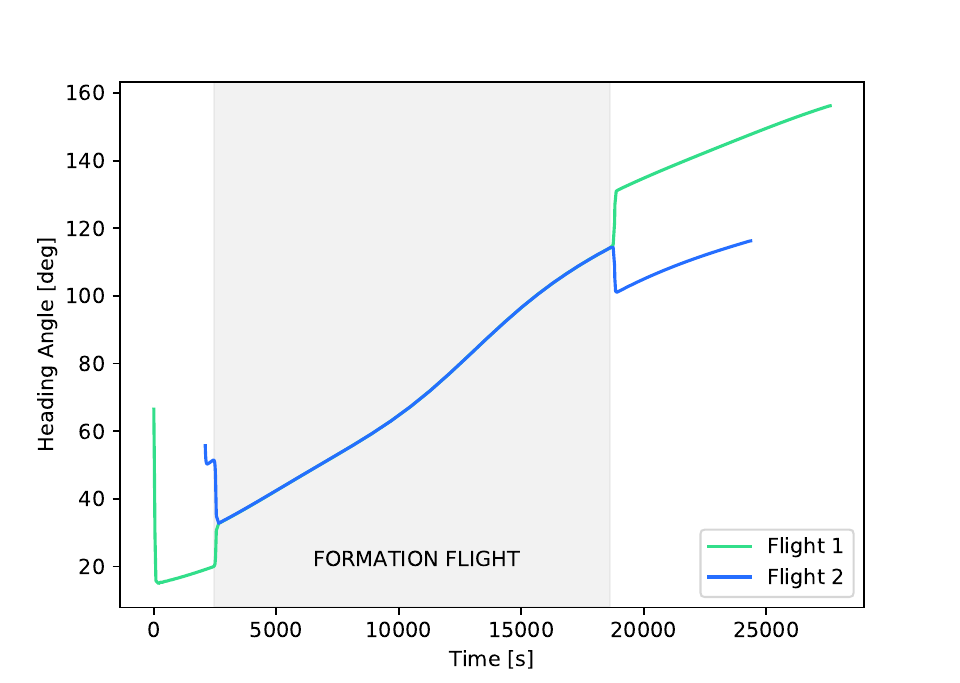}}
\hspace{-2mm}
\subfigure[Mach number]{\includegraphics[width=55mm]{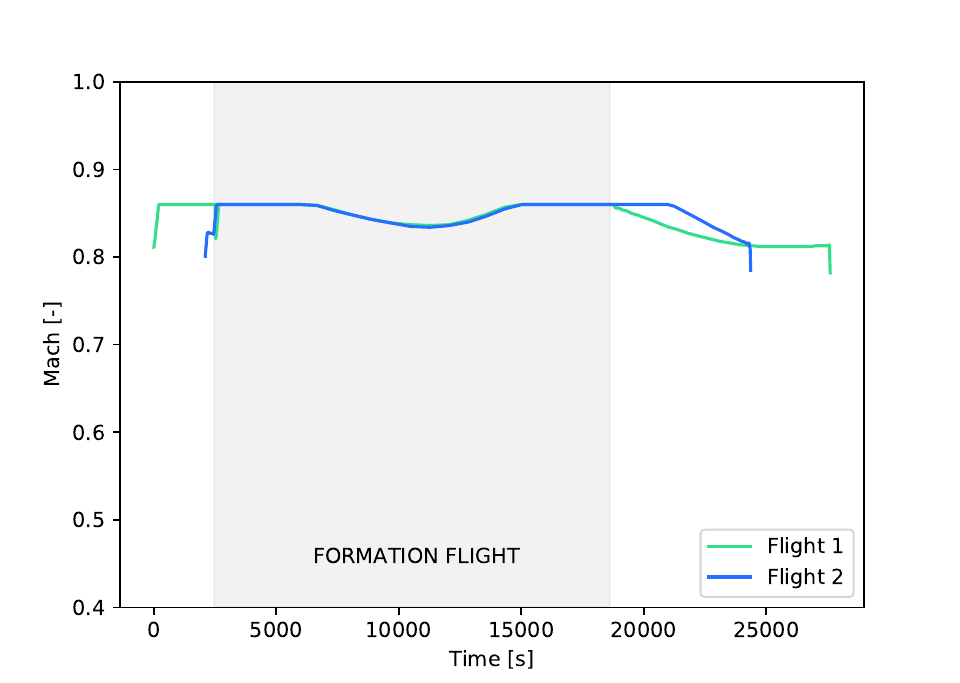}}
\caption{Experiment A: state variables.}
\label{fig:state_baseline}
\end{figure*}

\begin{figure*}[ht!]
\centering
\renewcommand{\figurename}{Fig.}
\hspace{-2mm}
\subfigure[Lift coeficient]{\includegraphics[width=55mm]{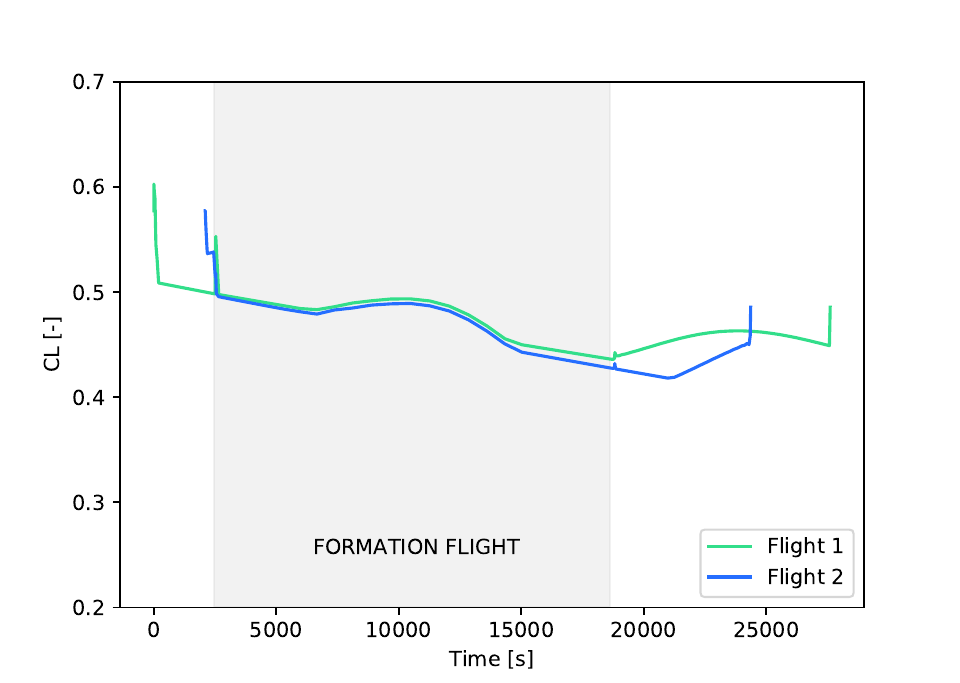}}
\hspace{-2mm}
\subfigure[Bank angle]{\includegraphics[width=55mm]{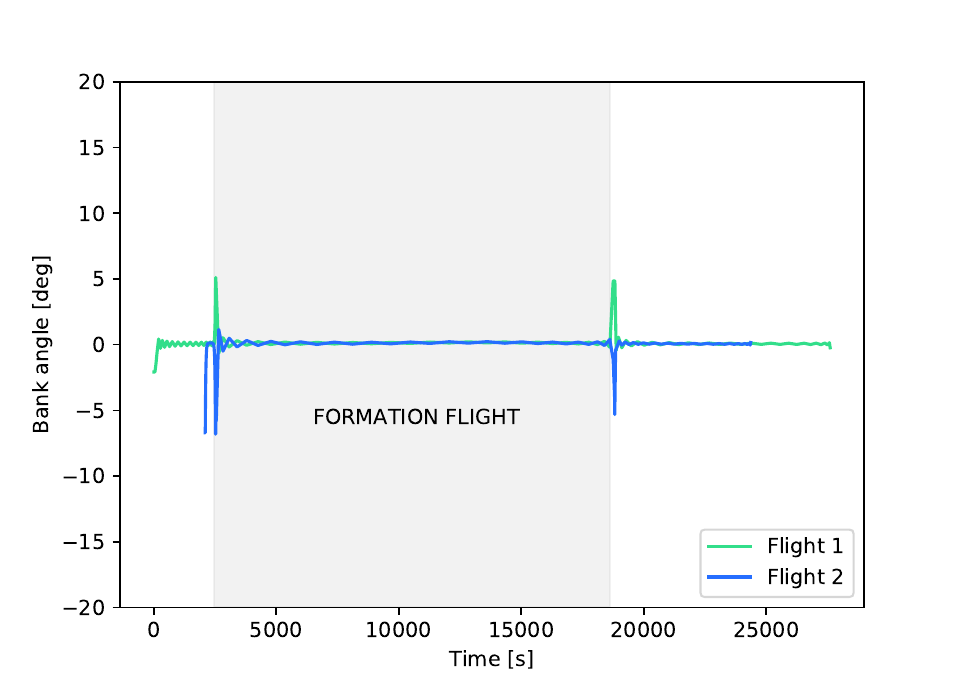}}
\subfigure[Adimensional thrust]{\includegraphics[width=55mm]{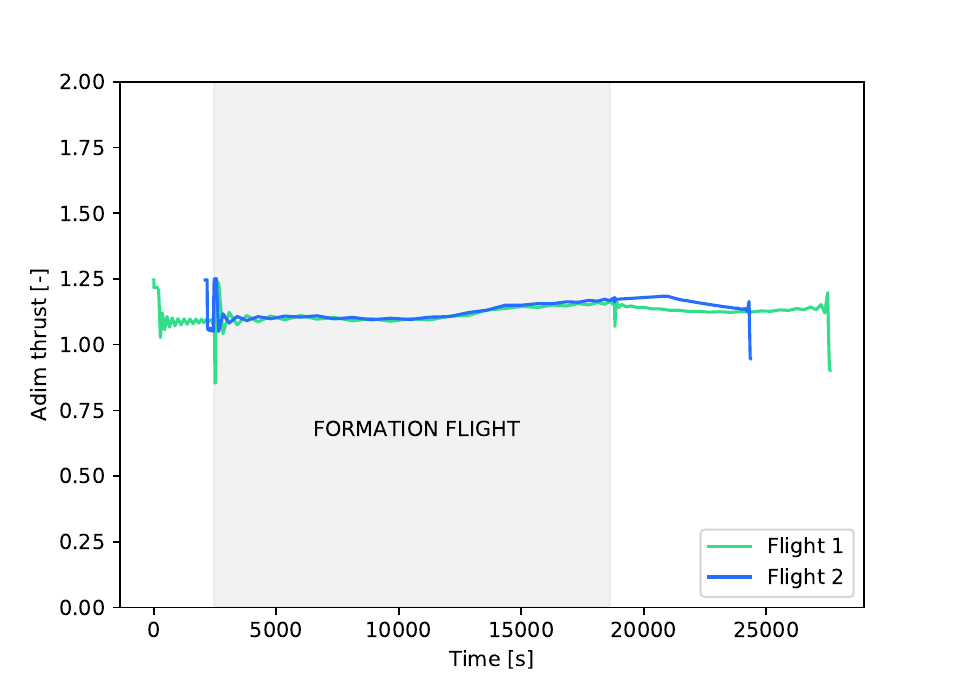}}
\caption{Experiment A: control variables.}
\label{fig:control_baseline}
\end{figure*}

The corresponding state and control variables are represented in Fig.~\ref{fig:state_baseline} and Fig.~\ref{fig:control_baseline}, respectively. For a better understanding of these figures, the portions of the plots of the variables that correspond to formation flight have been represented on a grey background.
It can be observed that the heading angle and the Mach number of both aircraft are the same during the formation, as opposed to the mass flow rate which is different for each aircraft. 
All the state variables vary smoothly, except the heading angle and the Mach number at the beginning and at the end of the flight. 
This behavior of the state variables is typical in boundary value problems, in which a quick maneuver is necessary to steer each aircraft from its initial state to its optimal route and from the optimal route to the final state. 
On the other hand, the control variables vary in a quite smooth way except 
around the switching instants. Notice that, the Gibbs phenomenon described in Section \ref{sect:solution_of_the_switched_optimal_control_problem} is present near the rendezvous and the splitting locations. However, this non-smooth behavior at the switching times lasts for a  very short period of time compared with the total flight time.

\begin{table}[ht!]
\centering
\caption{Experiment A: solo and formation flight results for Flight 1.}
\medskip
\resizebox{9cm}{!}
 {
\begin{tabular}{c cc cc cc c }
\multicolumn{1}{c }{}            & \multicolumn{1}{c }{\textbf{Flight} }  & \multicolumn{1}{c }{\textbf{Fuel}}     &\multicolumn{1}{c }{\textbf{Covered }  }  &\multicolumn{1}{c }{\textbf{DOC}  }   \\  
\multicolumn{1}{c }{}            & \multicolumn{1}{c }{\textbf{Time} \textbf{{[}}$\boldsymbol{h}$\textbf{{]}}}  & \multicolumn{1}{c }{\textbf{Burn} \textbf{{[}}$\boldsymbol{kg}$\textbf{{]}} }    &\multicolumn{1}{c }{\textbf{ Distance} \textbf{{[}}$\boldsymbol{km}$\textbf{{]}} }  &\multicolumn{1}{c }{ \textbf{{[}} \hspace{-0.2em}$\boldsymbol{mu}$\textbf{{]}}}   \\  
\cline{2-5}
\textbf{Formation F. }    &       7.68      &      46543.93   &   6716.44      & 40875.15  \\
\textbf{Solo F.}    &     7.47        &     48596.82    &     6205.81   & 42085.37   \\
  \hline
\end{tabular}
}
\label{table: results_comparison_baseline_solo_flight1}
\end{table}

\begin{table}[ht!]
\centering
\caption{Experiment A: solo and formation flight results for Flight 2.}
\medskip
\resizebox{9cm}{!}
 {
\begin{tabular}{c cc cc cc c }
\multicolumn{1}{c }{}            & \multicolumn{1}{c }{\textbf{Flight} }  & \multicolumn{1}{c }{\textbf{Fuel}}     &\multicolumn{1}{c }{\textbf{Covered }  }  &\multicolumn{1}{c }{\textbf{DOC}  }   \\  
\multicolumn{1}{c }{}            & \multicolumn{1}{c }{\textbf{Time} \textbf{{[}}$\boldsymbol{h}$\textbf{{]}}}  & \multicolumn{1}{c }{\textbf{Burn} \textbf{{[}}$\boldsymbol{kg}$\textbf{{]}} }    &\multicolumn{1}{c }{\textbf{ Distance} \textbf{{[}}$\boldsymbol{km}$\textbf{{]}} }  &\multicolumn{1}{c }{ \textbf{{[}} \hspace{-0.2em}$\boldsymbol{mu}$\textbf{{]}}}   \\  
\cline{2-5}
\textbf{Formation F. }    &      6.19    &   40130.77     &     5348.58     &  34776.74  \\
\textbf{Solo F.}    &      6.18       &     39683.14  &   5462.84  & 34452.60    \\
  \hline
\end{tabular}
}
\label{table: results_comparison_baseline_solo_flight2}
\end{table}

In this experiment, the great circle distance between departure and arrival locations of Flight 1 is 5761.08 km and the great circle distance between departure and arrival locations of Flight 2 is 5214.72 km. For these flights, the great circle distance between departure locations is 537.00 km and the great circle distance between arrival locations is 1244.77 km. Thus, this scenario and the selected fuel burn saving for the trailing aircraft can be considered not optimistic for formation flight.

Tables \ref{table: results_comparison_baseline_solo_flight1} and \ref{table: results_comparison_baseline_solo_flight2} summarize the results obtained in formation and solo flights for Flight 1 and Flight 2, respectively. As expected, the flight times increase for both aircraft in formation flight with respect to solo flight, although for Flight 2 this flight time increment is  negligible. Moreover, the fuel consumption of Aircraft 2, the leader aircraft, increases whereas the fuel consumption of Aircraft 1, the trailing one, decreases. The \textsf{DOC} is expressed in generic monetary units, {mu}. The total \textsf{DOC} of formation flight is almost one thousand generic monetary units lower than that of the solo flight, which amounts to more than one percent of reduction. 
For Flight 1, the total flight distance covered in formation flight is more than 500 km larger than in solo flight. Surprisingly, for Flight 2, the total flight distance covered in formation flight is smaller than the distance covered in solo flight. Obviously, both of them are higher than the great circle distance.
The reason behind this unexpected result for Flight 2 is the Jet Stream. Indeed, in order to get the greatest benefits from the wind field, Aircraft 2 in solo flight takes a route more than 200 km longer than the orthodromic route. In case of formation, this large detour, which would involve both aircraft, becomes less advantageous and therefore it does not appear in the optimal solution.

\subsection{Experiment B: two-aircraft transoceanic  mission design with delays in the departure times}

The flights considered in Experiment B, Flight 1 and Flight 2, have the same departure and arrival locations, the same boundary values for the state variables of each aircraft, and the same fuel burn savings for the trailing aircraft as in Experiment A. As before, in case of formation flight, Aircraft 1, associated to Flight 1, is forced to be the trailing aircraft and Aircraft 2, associated to Flight 2, is the leader one. The difference with respect to Experiment A is that several delays in the departure times of both flights have been considered.

This experiment has been carried out to determine how delays in the departure times of both flights affect the formation flight, in terms of routes, rendezvous and splitting locations and times, until formation flight becomes {non-optimal}  and, as a consequence, solo flights are selected by the algorithm.

\begin{table}[ht!]
\centering
\caption{Experiment B: delays in the departure time with respect to the baseline scenario for the different cases.}
\medskip
 {
\begin{tabular}{c c c}
\multicolumn{1}{c }{}            & \multicolumn{2}{c }{\textbf{Delays [min]} }   \\  
\cline{2-3}
\multicolumn{1}{ c }{\textbf{Cases}} & \textbf{Flight 1} & \textbf{Flight 2}    \\ \hline
$B_1$         &       20     &   0    \\
$B_2$         &       15     &   0        \\
$B_3$         &       10     &   0      \\
$B_4$         &       5      &   0      \\
\boldmath{$B_5$} &     \textbf{  - }     &   \textbf{-}      \\
$B_6$          &       0      &   5      \\
$B_{7}$       &       0      &   10      \\
$B_{8}$       &       0      &   15      \\
$B_{9}$       &       0      &   20      \\
$B_{10}$      &       0      &   25      \\
  \hline
\end{tabular}
}
\label{table:delays}
\end{table}

\begin{figure*}[ht!]
\centering
\renewcommand{\figurename}{Fig.}
\includegraphics[scale=0.5]{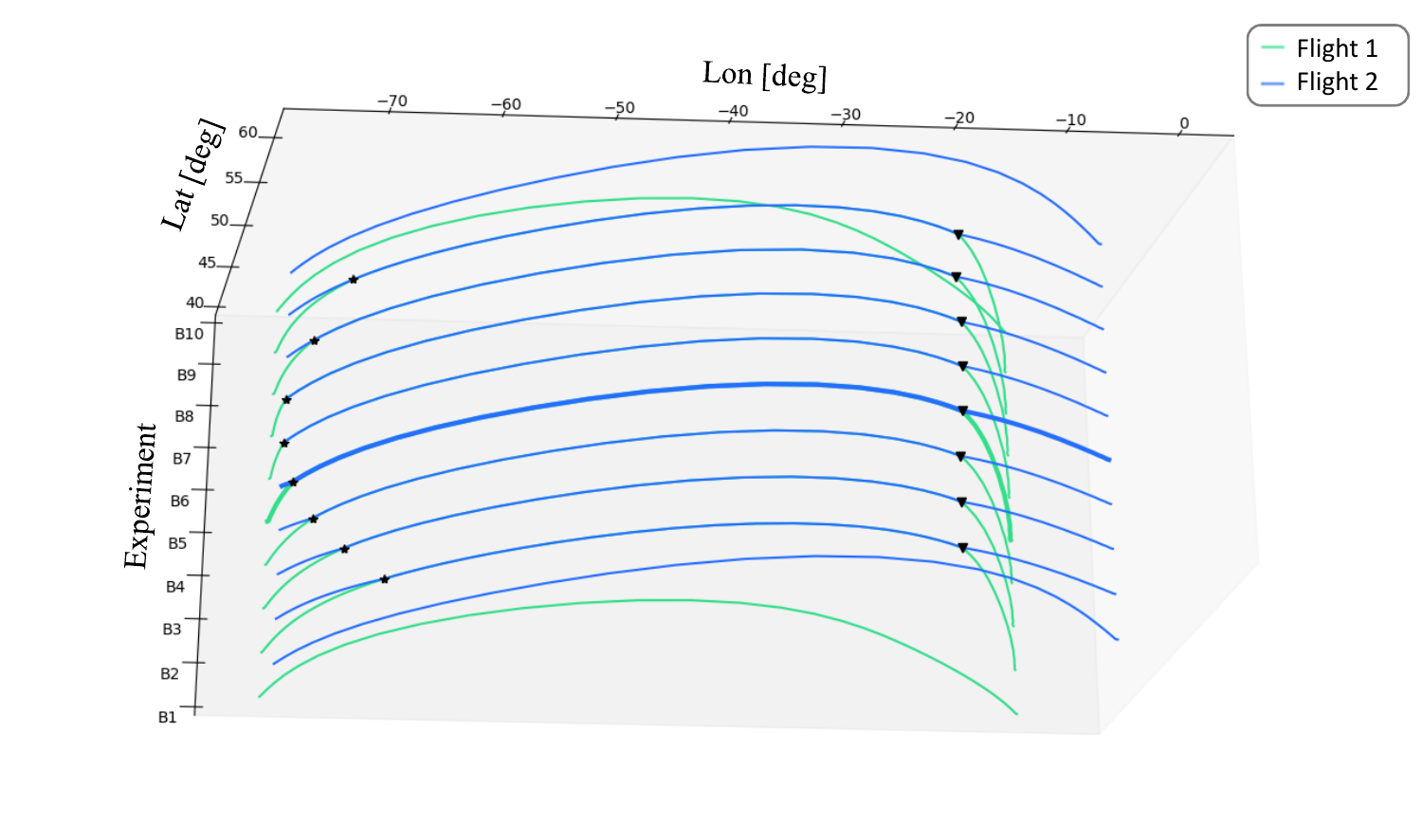}
\caption{Experiment B: routes obtained in the different cases considered.}
\label{fig:grafica_comparativa_tray_delays_RV_SP}
\end{figure*}

For this analysis, several numerical simulations have been conducted introducing delays in the departure times of both flights in Experiment A, namely, 10:15 for Flight 1 and 10:50 for Flight 2, which represents the baseline case. In particular, five-minute increments in the delays have been considered for each flight. The actual values of the delays are listed in Table \ref{table:delays}, in which each case is identified by a different symbol $B_1, \dots , B_{10}$, being $B_5$ the baseline case.

For a better understanding of the results, the optimal aircraft routes obtained for different delays in the departure times are represented together in Fig.~\ref{fig:grafica_comparativa_tray_delays_RV_SP}, in which the green and blue lines represent the optimal routes of Flight 1 and Flight 2, respectively. The optimal route obtained in the baseline case, $B_5$, is represented with a thicker line. Notice that, in this figure, the scale is not the same for latitude and longitude axes, and the symbols on the third axis denote the different cases.

In particular, Fig.~\ref{fig:grafica_comparativa_tray_delays_RV_SP} shows the rendezvous and splitting locations for each case, which are represented with small black stars and triangles, respectively. It can be observed that there is a great dependency of the location of the rendezvous points on the delay. For instance, in cases $B_6$ and $B_7$, the rendezvous location is very close to the departure location of Flight 2. On the contrary, the splitting locations barely change with delay.
When the delay on one departure time increases, since it takes more time for the delayed aircraft to reach the other one, the rendezvous point changes. However, once the formation is created, both aircraft follow a very similar route for all the cases regardless of the delay in the departure time, and therefore, the optimal splitting locations do not change significantly. In Table \ref{table:RV_SP_distances_ff}, the rendezvous and splitting locations and the distance covered flying in formation obtained in cases $B_2, \dots, B_9$ are displayed. Notice the proximity of the splitting locations in all the cases and the closeness among the rendezvous locations in cases $B_6$ and $B_7$ and the departure location of Flight 2.

On the contrary, the individual flight times and the duration of the formation flight, as well as the rendezvous and splitting times significantly change due to the initial detour made by both aircraft to create the formation. 
In Table \ref{table:RV_SP_times_ff}, the rendezvous and splitting times and the duration of the formation flight for each case are listed. 
Notice that in cases $B_1$ and $B_{10}$ solo flight is the optimal option. Therefore, the corresponding results have not been reported in tables and figures, where only the results obtained in formation flight are given. Further information on the results of cases $B_1$ and $B_{10}$ can be found in the results of Experiment A.

 \begin{table}[ht!]
\centering
\caption{Experiment B: rendezvous and splitting locations and formation flight distance.}
\medskip
\resizebox{9cm}{!}
 {
\begin{tabular}{cccc }
\multicolumn{1}{c }{}            & \multicolumn{1}{c }{\textbf{Rendezvous} }  & \multicolumn{1}{c }{\textbf{Splitting} }    &\multicolumn{1}{c }{\textbf{FF} }   \\
\multicolumn{1}{ c }{} & \textbf{(Lat, Lon)} \textbf{{[}}$\boldsymbol{deg}$\textbf{{]}}   & \textbf{(Lat, Lon)} \textbf{{[}}$\boldsymbol{deg}$\textbf{{]}}   &  \textbf{Distance} \textbf{{[}}$\boldsymbol{km}$\textbf{{]}}   \\ \hline
\boldmath{$B_2$}       &        (51.26, -65.69)    &     (57.05, -17.13)  &     3208.86   \\
\boldmath{$B_3$}       &        (49.15, -68.85)   &     (57.04,  -17.11)  &    3536.29   \\
\boldmath{$B_4$}       &       (47.10,-71.22)    &     (57.04,  -17.06)   &    3829.61   \\
\boldmath{$B_5$}       &        \textbf{(46.04, -72.86)}    &  \textbf{(57.05, -16.83)}     &    \textbf{4024.45}  \\
\boldmath{$B_6$}      &         (45.48, -73.71)    &     (57.06, -16.71)  &     4130.80   \\
\boldmath{$B_{7}$}    &           (45.49, -73.72) &     (57.06,  -16.71)   &     4130.68   \\
\boldmath{$B_{8}$}    &          (47.59, -71.91)  &     (57.06, -16.69)   &     3823.66   \\
\boldmath{$B_{9}$}    &         (49.92, -69.09)   &     (57.06, -16.73)  &     3523.89   \\
  \hline
\end{tabular}
}
\label{table:RV_SP_distances_ff}
\end{table}

\begin{table}[ht!]
\centering
\caption{Experiment B: rendezvous and splitting times and formation flight duration.}
\medskip
 {
\begin{tabular}{c c c c}
\multicolumn{1}{c }{}            & \multicolumn{3}{c }{\textbf{Time} \textbf{{[}}$\boldsymbol{h}$\textbf{{]}} }     \\  
\cline{2-4}
\multicolumn{1}{c }{}            & \multicolumn{1}{c }{\textbf{Rendezvous} }  & \multicolumn{1}{c } {\textbf{Splitting} }    &\multicolumn{1}{c }{\textbf{FF duration} }   \\ \hline
\boldmath{$B_2$}    &    1.58  &    5.21    &    3.63           \\
\boldmath{$B_3$}    &    1.22  &    5.21    &    3.99           \\
\boldmath{$B_4$}    &    0.90  &    5.22    &    4.33           \\
\boldmath{$B_5$}    &    \textbf{0.70} &   \textbf{5.24}    &   \textbf{4.54}  \\
\boldmath{$B_6$}    &    0.67  &    5.32    &   4.66            \\
\boldmath{$B_{7}$}  &    0.75  &    5.40    &   4.66            \\
\boldmath{$B_{8}$}  &    1.15  &    5.45    &   4.30            \\
\boldmath{$B_{9}$}  &    1.61  &    5.58    &   3.96            \\
  \hline
\end{tabular}
}
\label{table:RV_SP_times_ff}
\end{table}

\begin{table}[ht!]
\centering
\caption{Experiment B: flight times, fuel burn and \textsf{DOC} variations compared to solo flights.}
\medskip
 {
\begin{tabular}{c  cc  cc  c}
\multicolumn{1}{c }{}            & \multicolumn{2}{c }{\textbf{{$\boldsymbol{\Delta}$ \textbf{Time} \textbf{{[}}$\boldsymbol{\%}$\textbf{{]}}} } }  & \multicolumn{2}{c }{{$\boldsymbol{\Delta}$ \textbf{Fuel burn} \textbf{{[}}$\boldsymbol{\%}$\textbf{{]}}} }     & \multicolumn{1}{c }{   }  \\  
\cline{2-5}
\multicolumn{1}{ c }{} & \textbf{Flight 1} & \textbf{Flight 2}  & \textbf{Flight 1} & \textbf{Flight 2}   & {$\boldsymbol{\Delta}$ \textbf{DOC} \textbf{{[}}$\boldsymbol{\%}$\textbf{{]}}} \\ \hline
\boldmath{$B_2$}       &  +2.56     &     + 3.96       &    -3.30      &      +  1.97        &    -0.14    \\
\boldmath{$B_3$}       &  +2.56     &     + 2.60       &    -3.74      &      +  1.80        &    -0.51    \\
\boldmath{$B_4$}       &  +2.57     &     + 1.26       &    -4.13      &      +  1.55        &    -0.89     \\
\boldmath{$B_5$}     & \textbf{+2.79}    &    \textbf{+ 0.05}      &    \textbf{-4.22}      &      \textbf{+  1.13}    &    \textbf{-1.17}   \\
\boldmath{$B_6$}       &  +3.79     &     + 0.14       &    -4.15      &      +  1.01        &    -1.09      \\
\boldmath{$B_{7}$}     &  +4.94     &     + 0.13       &    -3.74      &      +  1.01        &    -0.79   \\
\boldmath{$B_{8}$}     &  +6.07     &     + 0.08       &    -3.06      &      +  1.01        &    -0.36   \\
\boldmath{$B_{9}$}     &  +7.25     &     + 0.07       &    -2.54      &      +  1.00        &    -0.02   \\
  \hline
\end{tabular}
}
\label{table:times_fuelburn_DOC_comparison}
\end{table}

The main results in terms of the objective functional are reported in Table \ref{table:times_fuelburn_DOC_comparison}. {In particular, the variations in the flight time and in the fuel consumption for each flight and the reduction in the \textsf{DOC} obtained in cases $B_2, \dots, B_9$, all of them compared to solo flights results, are reported. It can be observed that, in general, the more the flight time of Flight 1 increases, the less the flight time of Flight 2 does and vice versa.} 
Besides, for those cases in which Flight 2 is delayed, the flight time increment of Flight 2 is very small. About the fuel burn, it can be seen that the greatest benefits for the trailing aircraft are achieved in the baseline case, which is the optimal one, as established in Experiment A. Notice the similarity among increments in fuel burn for Flight 2 in the results of cases $B_6, B_7, B_8$, and $B_9$, in which Flight 2 has a delay in the departure time. The greatest \textsf{DOC} reduction is 1.17\% and the smallest is 0.02\%. This tiny \textsf{DOC} reduction implies that flying in formation in case $B_9$ has negligible benefits, which could be further reduced by any contingency during the flight.

In general, one may conclude that the greatest flight distances covered in formation flight
correspond to the highest total benefits in the fuel burn reduction. However, comparing the results displayed in tables \ref{table:RV_SP_distances_ff} and \ref{table:times_fuelburn_DOC_comparison}, it is easy to check that longer distances covered in formation flight and greater duration of formation flight do not always imply more benefits in terms of reduction of \textsf{DOC}. It can be seen in tables \ref{table:RV_SP_distances_ff} and \ref{table:RV_SP_times_ff}, that in the results of cases $B_6$ and $B_7$, aircraft are flying 4.66 hours and about 4130 km in formation, being the greatest duration of formation flights and the greatest distance covered in formation flights, respectively. Nevertheless, as it can be observed in Table \ref{table:times_fuelburn_DOC_comparison}, the corresponding reductions in the \textsf{DOC} with respect to solo flights are 1.09\% and 0.79\%, respectively. The reduction of the \textsf{DOC} is slightly smaller than the one obtained in the baseline case, which amounts to 1.17\%.

\begin{table}[ht!]
\centering
\caption{Experiment B: extra distance covered compared to solo flight distance.}
\medskip
\begin{tabular}{c  ccc}
\multicolumn{1}{c }{} & \multicolumn{2}{c }{{$\boldsymbol{\Delta}$ \textbf{Distance} \textbf{{[}}$\boldsymbol{km}$\textbf{{]}}} } \\   
\cline{2-3}
\multicolumn{1}{c }{}  & \multicolumn{1}{c }{Flight 1}  & \multicolumn{1}{c }{Flight 2}  \\   
\hline
   \boldmath{$B_2$}       & 451.21    & -123.32    \\ 
   \boldmath{$B_3$}       & 455.03    & -127.23    \\ 
   \boldmath{$B_4$}       & 467.05    & -123.25    \\ 
   \boldmath{$B_5$}       & \textbf{510.62}    & \textbf{-114.26}    \\ 
   \boldmath{$B_6$}       & 543.56    & -108.73    \\ 
   \boldmath{$B_7$}       & 544.79    & -111.76    \\ 
   \boldmath{$B_8$}       & 516.48    & -110.18    \\ 
   \boldmath{$B_9$}       & 507.60    & -109.73    \\ 
\hline
\end{tabular}
\label{table:extra_distance}
\end{table}

Table \ref{table:extra_distance} has been added to give information about the total detour done by each flight. As in Experiment A, in all the cases considered in Experiment B there is a reduction in the distance covered by Flight 2 in formation flight compared to solo flight. 
Noteworthy is the great detour made by Flight 1 in the formation flight solution: about 500 kilometers of diversion from the solo flight route.

\textcolor{black}{
The average computational time to find the solution has been $4.234$ s. For the sake of comparison, the same problem has been solved using a multiphase method \cite{hartjes2019trajectory}. In this case, two different \textsf{OCP} have been solved, namely, one in which aircraft flight solo and another one in which they are forced to fly in formation. The average computational time has been $7.784$ s.
}

\subsection{Experiment C: three-aircraft transoceanic mission design with different fuel savings schemes}

Experiment C involves three transoceanic eastbound flights, Flight 1, Flight 2 and Flight 3, with given fuel savings scheme and boundary values of the state variables. Flight 1, Flight 2, and Flight 3 are operated by Aircraft 1, Aircraft 2, and Aircraft 3, respectively.
The three flights considered in this experiment have the following departure and arrival locations, the first two of them being the same as in Experiment A:

\begin{itemize}
\item Flight 1: New York (\texttt{JFK}) - Madrid (\texttt{MAD}).
\item Flight 2: Montreal (\texttt{YUL}) - London (\texttt{LHR}).
\item Flight 3: Boston (\texttt{BOS}) - Paris (\texttt{CDG}).
\end{itemize}
The departure times of Flights 1,2, and 3, are set to 10:15, 10:50, and 10:30, respectively. The first two of them are the same as in Experiment A.
The boundary conditions for the state variables of the three aircraft are given in Table \ref{table:boundary_conditions_flight3}. The boundary conditions for the state variables of Aircraft 1 and  Aircraft 2 are the same as in Experiment A.

\begin{table}[ht!]
\centering
\caption{Experiment C: boundary conditions for the three flights.}
\medskip
\begin{tabular}{ccccc}
\multicolumn{1}{ c }{\textbf{Symbol}} & \textbf{{Units}}& \textbf{{Flight 1}}             & \textbf{{Flight 2}}  & \textbf{{Flight  3}}                   \\ \hline
 $\phi_I$ & [deg]    & 40.64  &  45.47   & 42.36 \\
 $\phi_F$ & [deg]    & 40.48  &  51.47   & 48.85  \\
 $\lambda_I$ & [deg] & -73.78 & -73.74   & -71.06  \\
 $\lambda_F$ & [deg] & -3.57  &  -0.45   & 2.35   \\
 $\chi_I$ & [deg]    & 66.51  &  55.70   & 56.46   \\
 $V_I$ & [m/s]       & 240    &   240   & 240   \\
 $V_F$ & [m/s]       & 220   &   220    & 220   \\
 $m_I$ & [kg]        & $220 \, 000$   &  $215 \, 000$  & $210 \, 000$   \\ \hline
\end{tabular}
\label{table:boundary_conditions_flight3}
\end{table}

Table \ref{table: results_third_solo_flight} summarizes the results obtained assuming that each flight is performed as a solo flight. It can be observed that the flight times, the fuel consumption, the covered distance, and the \textsf{DOC} of each flight are quite different.

\begin{table}[ht!]
\centering
\caption{Experiment C: results for the three solo flights.}
\medskip
\begin{tabular}{cccc}
\multicolumn{1}{l}{}  &\multicolumn{1}{c }{\textbf{Flight 1}}  & \multicolumn{1}{c }{\textbf{Flight 2}}  & \multicolumn{1}{c }{\textbf{Flight 3} }   \\  
\cline{1-4}
\textbf{Flight Time} \textbf{{[}}$\boldsymbol{h}$\textbf{{]}} &  7.47 &  6.19  &     6.90 \\
\textbf{Fuel burn} \textbf{{[}}$\boldsymbol{kg}$\textbf{{]}}  & 48596.82 &   39683.14 &   42877.40 \\
\textbf{Covered Distance} \textbf{{[}}$\boldsymbol{km}$\textbf{{]}}    &6205.81   & 5462.84 &      5945.89 \\
{\textbf{DOC} \textbf{{[}}$\boldsymbol{mu}$\textbf{{]}}}  &  42087.53   &  34462.32 &     37466.18     \\
\hline

\end{tabular}
\label{table: results_third_solo_flight}
\end{table}

An analysis is performed to determine how the fuel savings scheme affects the formation flight, in terms of routes, rendezvous and splitting locations and times, and flight times. 

As already mentioned, the only formation configuration allowed is the in-line formation and the relative position in the formation are fixed. The considered benefits for the intermediate and the trailing aircraft are ranging from 6\% to 14\%.
Several numerical simulations have been conducted introducing changes in the fuel savings for the intermediate and trailing aircraft, in which each case is identified by a different symbol $C_1,\dots,C_5$. The fuel savings considered in the different cases for the intermediate and trailing aircraft have been listed in Table \ref{table:fuel_savings_3aircraft}.

\begin{table}[ht!]
\centering
\caption{Experiment C: fuel savings for the intermediate and the trailing aircraft, in the different cases.}
\medskip
\begin{tabular}{c c }
\multicolumn{1}{c }{ }            & \multicolumn{1}{c }{\textbf{Fuel savings [\%]} }   \\  
 \hline
\boldmath{$C_1$}       &        6      \\
\boldmath{$C_2$}       &       8   \\
\boldmath{$C_3$}       &        10   \\
\boldmath{$C_4$}       &       12  \\
\boldmath{$C_5$}        &        14   \\
 \hline
\end{tabular}
\label{table:fuel_savings_3aircraft}
\end{table}

\begin{figure}[htb]
\centering
\renewcommand{\figurename}{Fig.}
\includegraphics[scale=0.48]{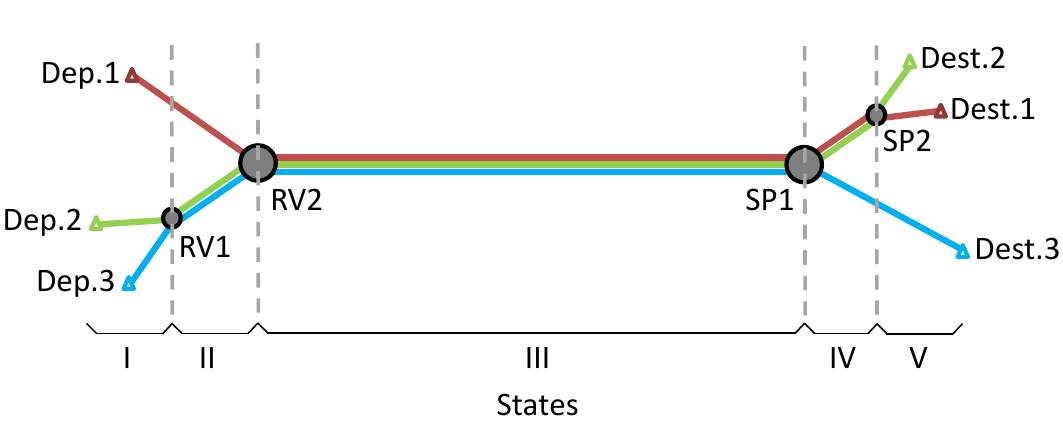}
\caption{Experiment C: discrete states representation for three-aircraft formation.}
\label{fig:five_discrete_states}
\end{figure}

It can be seen that in all the considered cases there are five discrete states of the switched dynamical system that represent the joint behavior of the aircraft. In Fig.~\ref{fig:five_discrete_states} a schematic representation of these five discrete states, from State I to State V, is given. In particular, each discrete state is characterized by the following modes

\begin{itemize}
\item[•] State I: all the aircraft fly in solo mode.
\item[•] State II: two aircraft fly in formation and one in solo mode.
\item[•] State III: all the aircraft fly in formation.
\item[•] State IV: two aircraft fly in formation and one in solo mode.
\item[•] State V: all the aircraft fly in solo mode.
\end{itemize}

It is easy to see that in the switches among these discrete states, there are two rendezvous points, \textsf{RV1} and \textsf{RV2}, where \textsf{RV1}  is the first one and \textsf{RV2}  is the second one, and two splitting points, \textsf{SP1}  and \textsf{SP2} , where \textsf{SP1}  is the first one and \textsf{SP2}  is the second one. In Fig.~\ref{fig:five_discrete_states}, the points \textsf{RV1}, \textsf{RV2}, \textsf{SP1}, and \textsf{SP2} are represented with grey points. Small triangles have been used to represent the departure and destination of each flight.

\begin{figure*}[ht!]
\centering
\renewcommand{\figurename}{Fig.}
\includegraphics[scale=0.5]{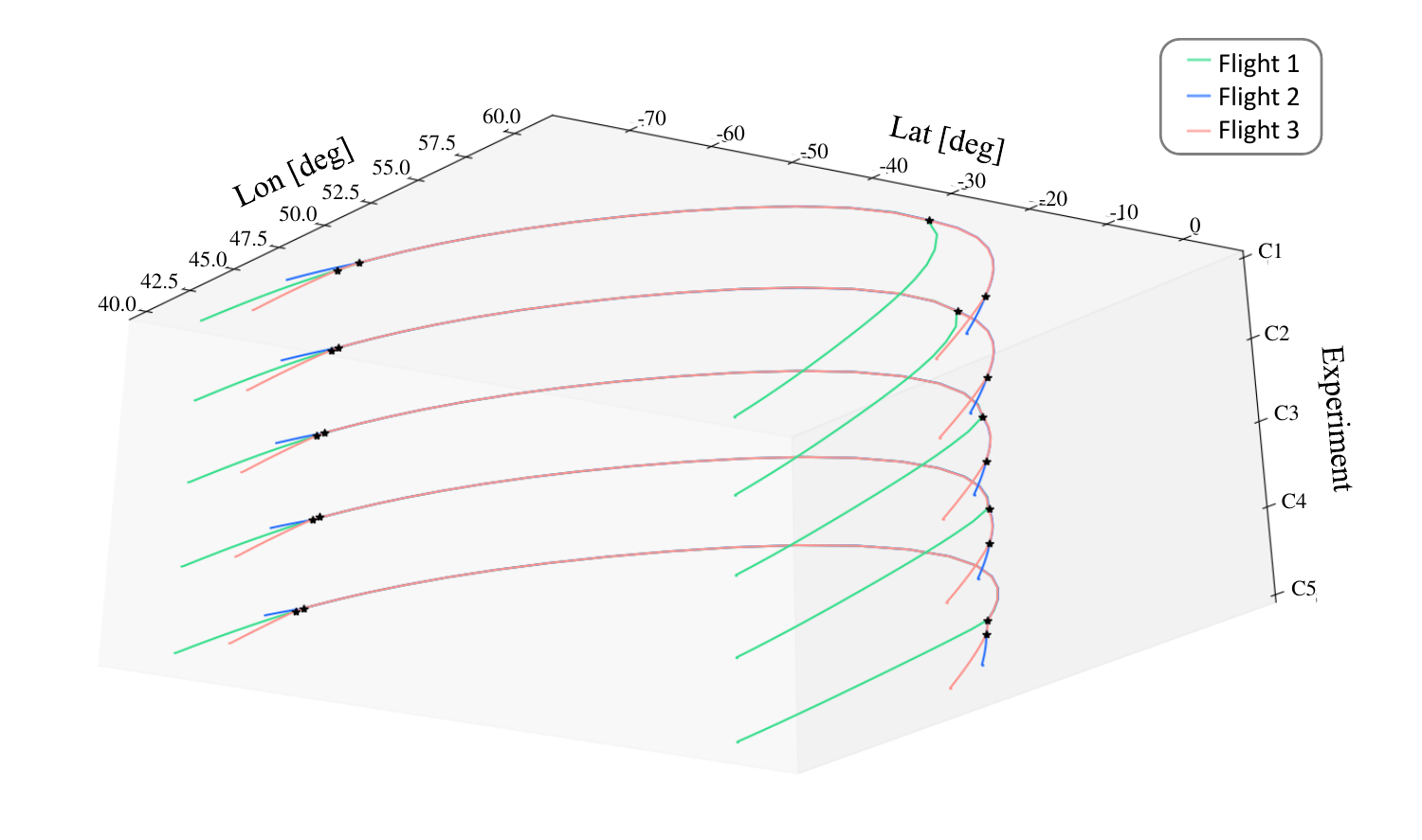}
\caption{Experiment C: routes obtained in the different cases considered. Detail of the rendezvous locations.}
\label{fig:exp3aircraft_RV}
\end{figure*}

\begin{figure*}[ht!]
\centering
\renewcommand{\figurename}{Fig.}
\includegraphics[scale=0.5]{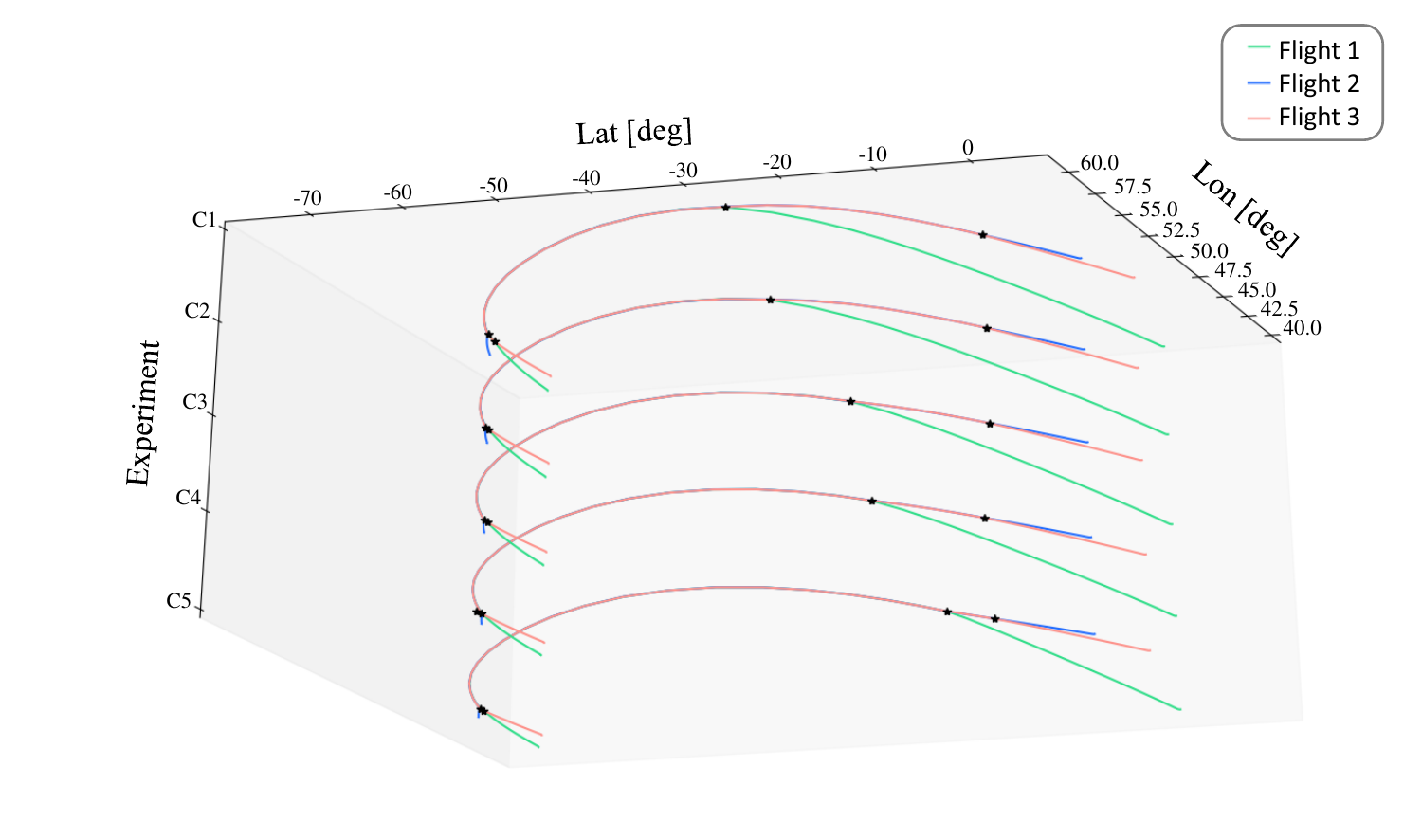}
\caption{Experiment C: routes obtained in the different cases considered. Detail of the splitting locations.}
\label{fig:exp3aircraft_SP}
\end{figure*}

Fig.~\ref{fig:exp3aircraft_RV} and Fig.~\ref{fig:exp3aircraft_SP} represent the three-aircraft routes obtained in the cases $C_1, \dots, C_5$. In these figures, the behavior of the rendezvous and splitting locations, respectively, {is} detailed. 
It can be observed that the location of both rendezvous points, \textsf{RV1} and \textsf{RV2}, and the location of the second splitting point, \textsf{SP2}, do not present significant variations in the different cases. On the contrary, the location of the first splitting point, \textsf{SP1}, notably changes, differing in some cases in more than 1000 km. This variation in the location of \textsf{SP1}, which is the point in which there is a switch between State III and State IV,  implies that the three-aircraft formation time and distance largely depend on the fuel savings scheme. Consequently, the formation time and distance in State IV also have a high dependency on it. Results in tables \ref{table:ff_distances_3aircraft} and \ref{table:ff_times_3aircraft} confirm this conclusion. 
It can also be observed in these tables that, in all the considered cases, the formation times and distances corresponding to discrete State II are small compared to other discrete states.

\begin{table}[ht!]
\centering
\caption{Experiment C: formation distances for the two- and three-aircraft formation phases.}
\medskip
\begin{tabular}{c c c c }
\multicolumn{1}{c }{}            & \multicolumn{3}{c }{\textbf{Formation distances [$\boldsymbol{km}$]} }    \\  
\cline{2-4}
\multicolumn{1}{ c }{\textbf{Cases}}   &  \multicolumn{1}{ c }{\textbf{State II}} &  \multicolumn{1}{ c }{\textbf{State III}} &  \multicolumn{1}{ c }{\textbf{State IV}} \\ \hline
\boldmath{$C_1$}       &  113.53    & 3008.99       &   1371.24   \\
\boldmath{$C_2$}       &   38.19    & 3323.62       &   1142.36   \\
\boldmath{$C_3$}       &   41.01    & 3743.05       &    737.69       \\
\boldmath{$C_4$}      &   44.48    & 3847.74       &    595.32   \\
\boldmath{$C_5$}       &   41.44    & 4229.41       &    257.29     \\
  \hline
\end{tabular}
\label{table:ff_distances_3aircraft}
\end{table}

\begin{table}[ht!]
\centering
\caption{Experiment C: formation times for the two- and three-aircraft formation phases.}
\medskip
\begin{tabular}{c c c c }
\multicolumn{1}{c }{}            & \multicolumn{3}{c }{\textbf{Formation times [$\boldsymbol{h}$]} }    \\  
\cline{2-4}
\multicolumn{1}{ c }{\textbf{Cases}}   &  \multicolumn{1}{ c }{\textbf{State II}} &  \multicolumn{1}{ c }{\textbf{State III}} &  \multicolumn{1}{ c }{\textbf{State IV}} \\ \hline
\boldmath{$C_1$}       &  0.13   & 3.30      &  1.72   \\
\boldmath{$C_2$}       &  0.06   & 3.69      &  1.43   \\
\boldmath{$C_3$}       &  0.05   & 4.26      &  0.92       \\
\boldmath{$C_4$}       &  0.04   & 4.42      &  0.84   \\
\boldmath{$C_5$}       &  0.06   & 4.92      &  0.74     \\
  \hline
\end{tabular}
\label{table:ff_times_3aircraft}
\end{table}

\begin{table*}[hbt]
\centering
\caption{Experiment C: total flight times, fuel burn and \textsf{DOC} reduction compared to solo flights, in the different cases.}
\medskip
 {
\begin{tabular}{c  ccc ccc c}
\multicolumn{1}{c }{}            & \multicolumn{3}{c }{\textbf{{$\boldsymbol{\Delta}$ \hspace{-0.4em} \textbf{Time} \textbf{{[}}$\boldsymbol{\%}$\textbf{{]}}} } }   & \multicolumn{3}{c }{\textbf{{$\boldsymbol{\Delta}$ \hspace{-0.4em} \textbf{Fuel burn} \textbf{{[}}$\boldsymbol{\%}$\textbf{{]}}} } }  & \multicolumn{1}{c }{}            \\  
\cline{2-8}
\multicolumn{1}{ c }{ } & \textbf{Flight 1} & \textbf{Flight 2}  & \textbf{Flight 3} & \textbf{Flight 1} & \textbf{Flight 2}  & \textbf{Flight 3}  & \multicolumn{1}{c }{\textbf{{$\boldsymbol{\Delta}$\hspace{-0.3em} \textbf{DOC} \textbf{{[}}$\boldsymbol{\%}$\textbf{{]}}} }} \\ \hline 
\boldmath{$C_1$}      &   2.28    &    0.13  &   0.30   &  -1.83   &    -4.17  &  0.90   & -1.14     \\
\boldmath{$C_2$}       &   2.45    &   -0.02  &   0.16   &  -2.87   &    -5.81  &  1.10  & -1.80     \\
\boldmath{$C_3$}       &   2.89    &   -0.06  &   0.10   &  -4.15   &    -7.39  &  1.25  & -2.51      \\
\boldmath{$C_4$}       &   2.90    &   -0.19  &   0.01  &  -5.29   &    -9.14  &  1.41     & -3.25  \\  
\boldmath{$C_5$}       &   3.54    &   -0.10  &   0.02   &  -6.69   &   -10.64  &  1.45   & -3.97   \\
\hline
\end{tabular}
}
\label{table:times_fuel_comparison_3aircraft}
\end{table*}

The main results related to the objective function obtained in Experiment C are reported in Table \ref{table:times_fuel_comparison_3aircraft}. In particular, the time increment of each flight, the reduction or increase of the fuel burn and the reduction in the \textsf{DOC} with respect to the results obtained assuming each flight is performed as a solo flight, are reported.  \textcolor{black}{It} can be seen in this table that, for case $C_5$, the reduction in the \textsf{DOC} is almost 4\% comparing to the case in which formation is not allowed.

\textcolor{black}{
The average computational time to find the solution has been $15.274$ s. For the sake of comparison, the same problem has been solved using a multiphase method \cite{hartjes2019trajectory}. In this case, 13 different \textsf{OCP} have been solved, since, in this case, there are 13 possible flight phase sequencing options, as shown in \cite{hartjes2019trajectory}.  
The average computational time has been $90.948$ s. 
This shows that the method presented in this paper is able to drastically reduce the computational time to solve formation mission design problems. 
}

\section{Conclusions}
\label{sect:conclusions}

In this paper, 
a novel framework to solve the formation mission design problem is presented. 
In the proposed framework, the mission is modeled as a switched dynamical system, in which aircraft are assumed to have two flight modes, namely solo and formation flight, and the discrete state of the switched dynamical system  is the result of their combination.
The discrete dynamics is modeled using logical constraints in disjunctive form. The formation mission design problem is transformed into an optimal control problem for a switched dynamical system with logical constraints, which is solved using the embedding approach.

The embedding approach is a unifying technique able to efficiently tackle both the switching dynamics and the logical constrains of the optimal control problem, transforming it into a smooth optimal control problem, which has been solved using a knotting pseudospectral method.

The main advantages of this approach are that the multi-phase formulation is avoided, as well as the use of binary
variables, decreasing the computational time and effort in finding the solution. This approach is substantially different from the previous approaches, which are based on exhaustively analyzing every possible formation mission individually and then, comparing the results.
Additionally, this approach is easily scalable to design formation mission problems
of an arbitrary number of aircraft.

Several numerical experiments have been conducted with two and three transoceanic flights. The results demonstrate that the proposed framework for mission design is fast and accurate and that the formation flight has great potential to reduce fuel consumption and emissions and, therefore, to mitigate the environmental impact of air transport sector. An analysis of the solutions has also been carried out in order to study how the delays in the departure time and changes in the fuel saving scheme influence the formation mission. The results indicate that delays and changes in the fuel savings scheme have significant influence on the formation mission.

This suggests that a deeper understanding of the influence of uncertainties in these and other parameters, such as the initial mass of the aircraft or meteorological conditions, on the formation mission is needed to assess the feasibility of formation flight in realistic scenarios usually characterized by in the presence of uncertainties. This will be subject of future research.



\section{Acknowledgments} \label{sect:acknowledgments}
This work has been partially supported by the grants number \texttt{TRA2017-91203} \texttt{-EXP} and \texttt{RTI2018-098471-B-C33} of the Spanish Government.

\bibliographystyle{IEEEtran}
\bibliography{IEEEabrv,biblio_maria}

\begin{IEEEbiography}
[{\includegraphics[width=1in,height=1.5in,clip,keepaspectratio]{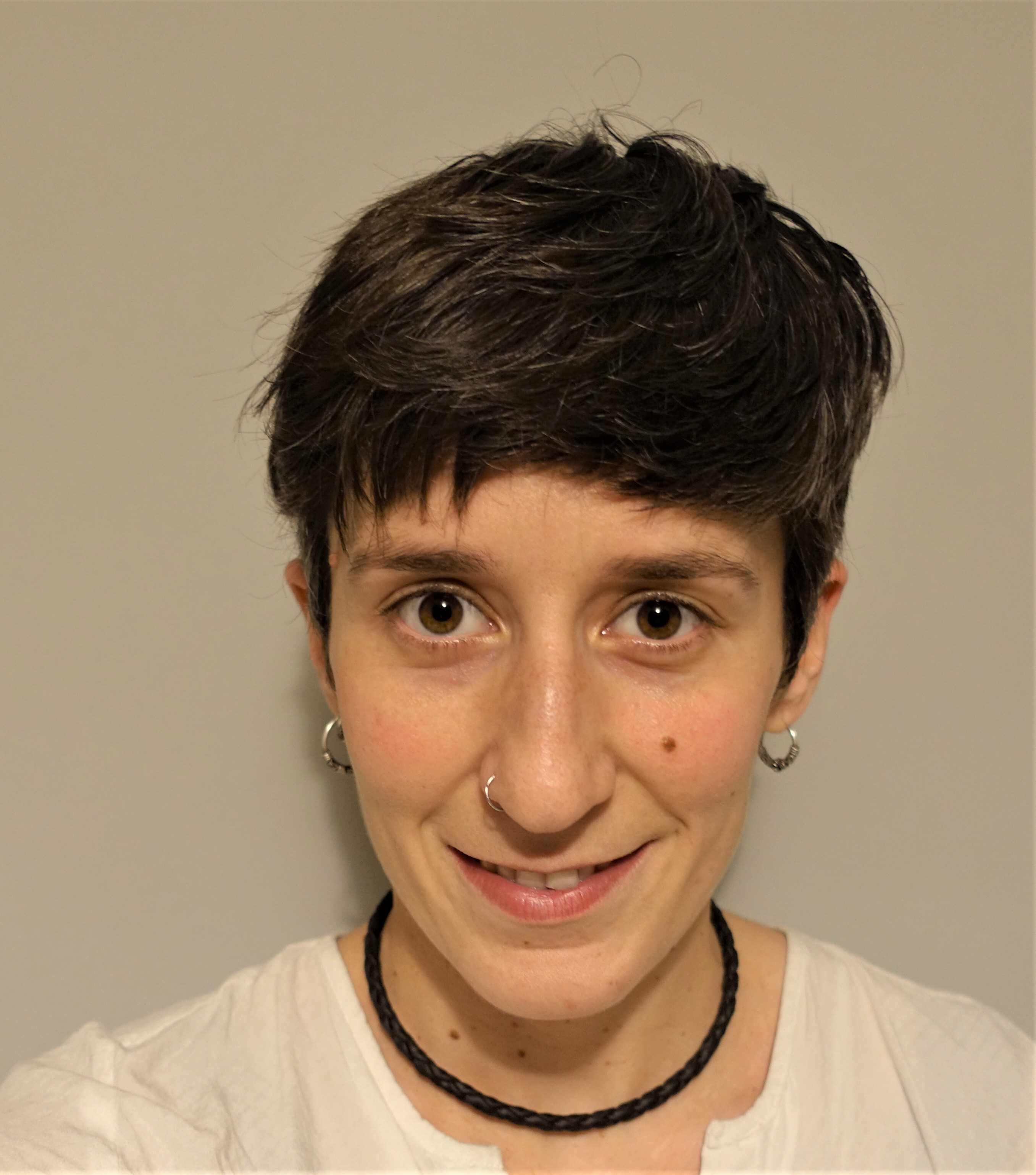}}]
{María Cerezo-Magaña} is a Teaching Assistant of {Aerospace Engineering} and a PhD student at the Universidad Rey Juan Carlos in Madrid, Spain. She received her MSc degree in Aeronautical Engineering from the Universidad Polit\'ecnica de Madrid. Her research is focused on deterministic and stochastic hybrid optimal control applied to aircraft trajectory optimization.
\end{IEEEbiography}

\begin{IEEEbiography}
[{\includegraphics[width=1in,height=1.25in,clip,keepaspectratio]{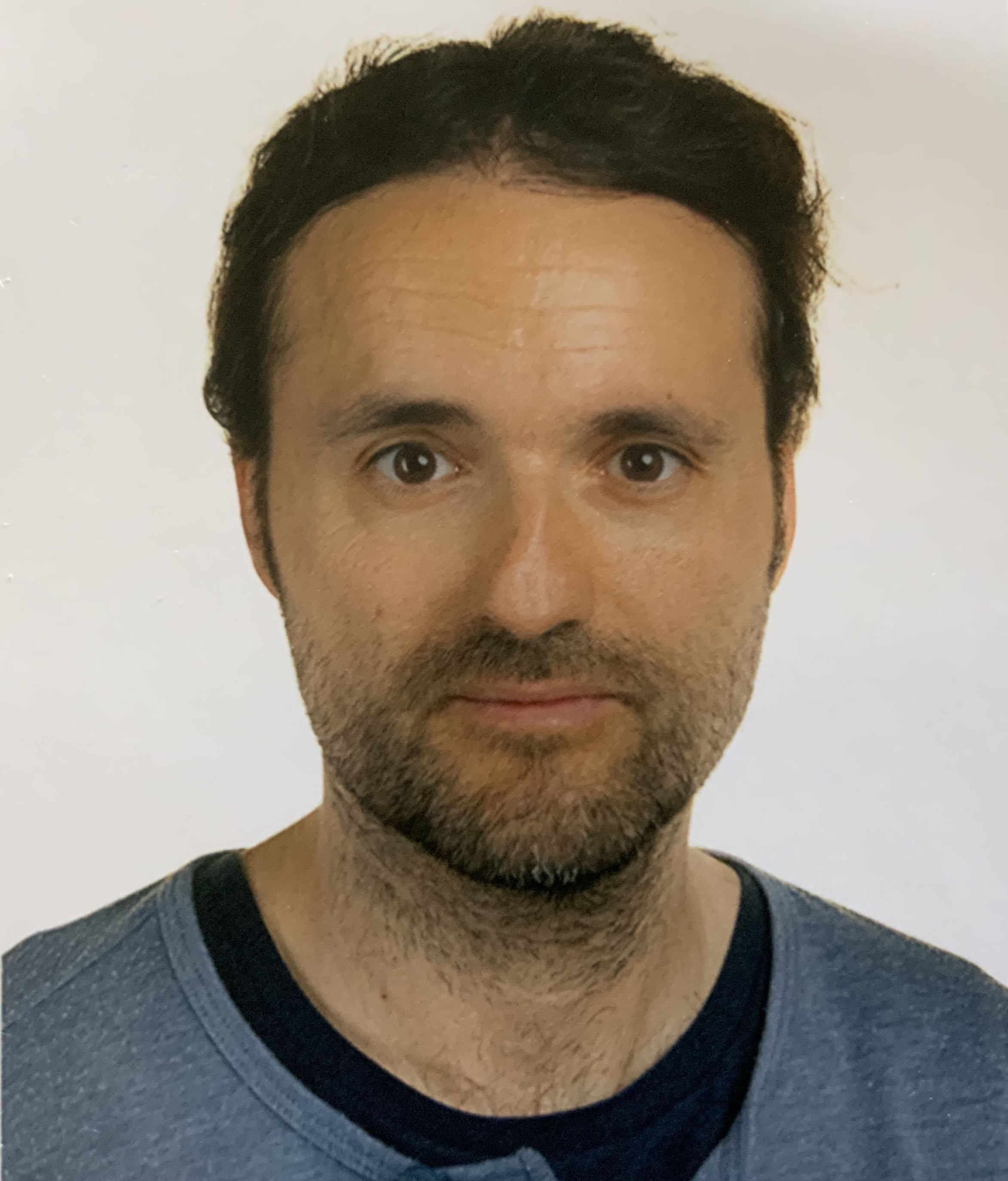}}]{Alberto Olivares}
is a Professor of Statistics and Vector Calculus at the Universidad Rey Juan Carlos in Madrid, Spain. He
received his MSc degree in Mathematics and his BSc degree in Statistics from the Universidad de Salamanca, Spain, and his PhD degree in Mathematical Engineering from the Universidad Rey Juan Carlos. He worked with the Athens University of Economics and Business. His research interests include statistical learning, stochastic hybrid optimal control and model predictive control with applications to biomedicine, robotics, aeronautics and astronautics.
\end{IEEEbiography}

\begin{IEEEbiography}
[{\includegraphics[width=1in,height=1.25in,clip,keepaspectratio]{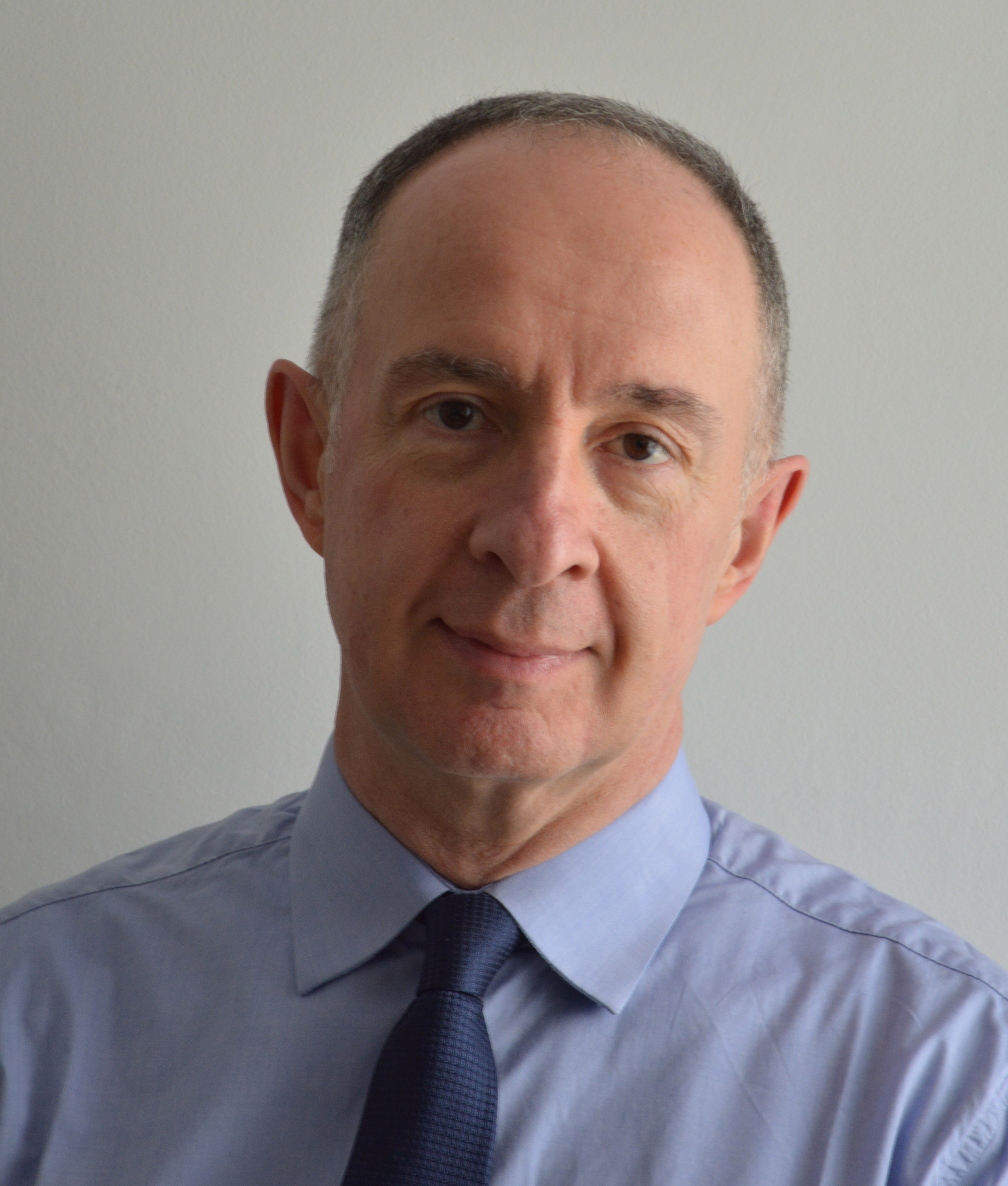}}]{Ernesto Staffetti}
is a Professor of Statistics and Control Systems at the Universidad Rey Juan Carlos in Madrid, Spain. He
received his MSc degree in Automation Engineering from the Universit\`a degli Studi di Roma ``La Sapienza,''
and his PhD degree in Advanced Automation Engineering from the Universitat Polit\`ecnica de Catalunya. 
He worked with the Universitat Polit\`ecnica de Catalunya, the
Katholieke Universiteit Leuven, the Spanish Consejo Superior de
Investigaciones Cient\'{\i}ficas, and with the University of North Carolina
at Charlotte. His research interests include stochastic hybrid optimal
control, iterative learning control and model predictive control with applications to robotics, aeronautics and astronautics.
\end{IEEEbiography}

\vfill

\end{document}